%% file: svungoc.tex
\documentclass[a4paper,11pt]{article}

\usepackage{ifthen}%
\newboolean{PDF}%
\setboolean{PDF}{false} 

\usepackage[latin1]{inputenc}
\usepackage{epsfig}
\usepackage{amsmath,amsfonts,amssymb}

\usepackage{mathptmx} 

\input{macro}
\newcommand{\C}{\mathcal{C}}
\renewcommand{\A}{\mathbb{A}}
\newcommand{\seq}{\stackrel{\mathrm{s}}{\sim}}
\newcommand{\h}{\hbar}
\renewcommand{\oh}{\mathcal{O}(\h^\infty)}
\newcommand{\pdos}{pseudodifferential operators}
\newcommand{\aangles}{Liou\-ville-Arnold-Mi\-neur}

\title{Symplectic Techniques for Semiclassical Completely Integrable
  Systems}

\author{V{\small\~U} N{\small G\d OC} San}
\date{}

\begin{document}

\maketitle

{\small\begin{center}
    Institut Fourier, Universit\'e Grenoble 1\\BP 74\\38402 Saint
    Martin d'H\`eres Cedex, France\\
    \texttt{san.vu-ngoc@ujf-grenoble.fr}\\
    \verb|http://www-fourier.ujf-grenoble.fr/~svungoc|
\end{center}}

\begin{abstract}
  This article is a survey of classical and quantum completely
  integrable systems from the viewpoint of local ``phase space''
  analysis. It advocates the use of normal forms and shows how to get
  global information from glueing local pieces. Many crucial phenomena
  such as monodromy or eigenvalue concentration are shown to arise
  from the presence of non-degenerate critical points.
\end{abstract}

\section{Foreword}

This article is mainly an adaptation and a translation of the last
chapters of my habilitation thesis defended in December 2003. Its aim
is to describe, in a unified way, old and new results in the theory of
classical and quantum (or rather semiclassical) completely integrable
systems in the spirit of the famous Darboux-Carathéodory theorem. This
theorem (which was essentially already known to Liouville) gives a
\emph{symplectic local normal form} for a classical completely
integrable system near a regular point. From the viewpoint of modern
geometers, it is very natural to build up the global theory from such
local results. It was far less obvious to apply this idea to quantum
systems, until the appearance of \pdos{}, \fio{}s and
microlocalisation techniques in the late 1960's. Nowadays, if one is
reasonably familiar with both geometric and microlocal techniques, it
seems evident that, as we did for classical systems, one should be
able to discover the global theory of semiclassical integrable systems
from local semiclassical analogues of the Darboux-Carathéodory
theorem. This, in essence, is what this article is about. Although it
is not my purpose to give the reader many details (especially about
semiclassical theories), I still hope that the text manages to convey
the right intuition. A more complete treatment will be published
elsewhere~\cite{san-panoramas}.

\section{What is a completely integrable system ?}

While the standard notion of completely integrable systems is now
perfectly standard, it is not obvious at all to define the
\emph{moduli space} of all integrable systems. More simply put: when
shall we say that two completely integrable systems are equivalent ?

\subsection{Classical mechanics}

Let $M$ be a $\Cinf$ symplectic manifold of dimension $2n$. The
algebra $\Cinf(M)$ of classical observables (or Hamiltonians) is
equipped with the symplectic Poisson bracket $\{\cdot,\cdot\}$. Any
function $H\in\Cinf(M)$ gives rise to a vector field $\ham{f}$ (the
Hamiltonian vector field) which, as a derivation, is just the Poisson
bracket by $H$; in other words the evolution of a function $f$ under
the flow of $\ham{H}$ is given by the equation
\[
\dot{f}=\{H,f\} .
\]
An \emph{integral} of the Hamiltonian $H$ is a function which is
invariant under the flow of $\ham{H}$; this means a function $f$ such
that $\{H,f\}=0$. The Hamiltonian $H$ is called \emph{completely
  integrable} if there exists $n-1$ independent functions
$f_2,\dots,f_n$ which are integrals of $H$ ($\{H,f_j\}=0$) and
moreover \emph{pairwise commute}: $\{f_i,f_j\}=0$. (This last
condition is always a consequence of the former for ``generic'' $H$;
but since we are going to study particular models or \emph{normal
  forms}, which are not generic, this condition is crucial.)

Actually, one sees from the definition that the function $H$ does not
play any distinguished role among the other functions $f_2,\dots,f_n$.
Our point of view will always be to consider, as a whole, a
\emph{classical completely integrable system} to be the data of $n$
functions $f_1,\dots,f_n$ in involution: $\{f_i,f_j\}=0$, which are
independent in the sense that for almost every point $m\in M$,
$df_1(m),\dots,df_n(m)$ are linearly independent. $n$ is the largest
number of such functions for which this is indeed possible.

We define the \emph{momentum map} of the system to be the map
$F=(f_1,\dots,f_n):M\fleche\RM^n$. In the language of Hamiltonian Lie
group actions, this is indeed a momentum (or ``moment'') map for a
local action of $\RM^n$ on $M$.

It is tempting to say that two completely integrable systems are
equivalent when their momentum maps are equivalent, in the sense that
there exists a diffeomorphism $g$ of $\RM^n$ such that $F_1=g\circ
F_2$. However, as we shall see later, this notion is too strong,
especially due to the existence of flat functions in the $\Cinf$
category.

The first natural attempt to weaken this equivalence is the following.
\begin{defi}
  Let $U$ be an open subset of $M$.  The \textbf{momentum algebra}
  $\mathbf{f}_{\restr U}=\gener{f_1,\dots,f_n}$ is the linear span of
  the $f_j$'s, as an abelian subalgebra of $\Cinf(U)$. The
  \textbf{commutant} of $\mathbf{f}$ is the set of all $g\in\Cinf(U)$
  Poisson-commuting with $\mathbf{f}$. It is denoted by
  $\C_\mathbf{f}(U)$~.
\end{defi}
Notice that by Jacobi's identity, $\C_\mathbf{f}(U)$ is a Lie algebra
under Poisson bracket.
\begin{defi}
  \label{defi:eq-faible}
  Let $F=(f_1,\dots,f_n)$ and $G=(g_1,\dots,g_n)$ be two completely
  integrable systems with momentum algebras $\mathbf{f}$ and
  $\mathbf{g}$. Then $F$ and $G$ are called \textbf{weakly equivalent}
  on an open subset $U$ if
  \[
  \C_{\mathbf{f}}(U)=\C_{\mathbf{g}}(U).
  \]
\end{defi}
It is clear that the definition of weak equivalence does not depend on
the particular choice of basis for $\mathbf{f}$ and $\mathbf{g}$.
Therefore, by a slight abuse of notation, the fact that the two
systems in consideration are equivalent shall be denoted as
\[
\mathbf{f}\sim_U\mathbf{g} .
\]
The associations $U\fleche\mathbf{f}_{\restr U}$ or $U\fleche
\C_\mathbf{f}(U)$ define typical presheaves over $M$. With this in
mind, we shall often simply refer to $\mathbf{f}$ or $\C_\mathbf{f}$
when the localisation on $U$ is unimportant or clearly implicit.  See
also corollary~\ref{coro:restriction}.  

The geometric interpretation of this definition is clear: The fibres
of $F$ define on $M$ a \emph{singular Lagrangian foliation}: when
$c\in\RM^n$ is a regular value of $F$ ($\textup{rank }dF(m)=n$ for all
$m\in F^{-1}(c)$), then $F^{-1}(c)$ is a Lagrangian submanifold of
$M$. The leaves of the singular Lagrangian foliation are all the
connected components of the fibres of $F$. Then $\C_\mathbf{f}$ is
just the algebra of smooth functions that are constant on the leaves.
So $\mathbf{f}\sim\mathbf{g}$ says that the spaces of leaves of the
corresponding singular foliations are the same for $\mathbf{f}$ and
$\mathbf{g}$, in a smooth way.

It turns out that the weak equivalence is not able, in some cases, to
distinguish singularities of integrable systems; this is due to the
fact that there is no requirement that the functions $f_1,\dots,f_n$
always be a \emph{reduced} set of equations for the foliation (in the
algebraic geometry sense). For instance the functions $x$ and $x^3$ on
$\RM^2$ give weakly equivalent systems, while as functions they
obviously don't have the same singularity type.

\begin{defi}
  \label{defi:eq-forte}
  Let $F=(f_1,\dots,f_n)$ and $G=(g_1,\dots,g_n)$ be two completely
  integrable systems with momentum algebras $\mathbf{f}$ and
  $\mathbf{g}$. Let $m\in M$. Then $F$ and $G$ are called
  \textbf{strongly equivalent} at $m$ if there exists a small
  neighbourhood of $m$ on which
  \[
  \C_{\mathbf{f}}.(\mathbf{f}-\mathbf{f}(m))=
  \C_{\mathbf{g}}.(\mathbf{g}-\mathbf{g}(m)).
  \]
  This will be denoted as $\mathbf{f}\seq\mathbf{g}$~.
\end{defi}
It is elementary (and probably standard in algebraic geometry) to see
that two strongly equivalent systems have exactly the same singularity
type. To do this, given a momentum algebra $\mathbf{f}$ vanishing at
$m$, we shall say that a basis $(f_1,\dots,f_n)$ of $\mathbf{f}$ if
well ordered when there is a partition of $n=n_1+\cdots+n_d$, $n_k\geq
0$ such that for all $k\in [1..d]$, the vector space spanned by the
$f_j$ corresponding to the block $n_k$ (which means
$j\in[(n_1+\cdots+n_{k-1}+1)..(n_1+\cdots+n_k)$) consists of functions
vanishing exactly at order $k$ at $m$ (except for the zero function).
The associated partition is unique. For instance if $m$ is a regular
point, every basis of $\mathbf{f}$ is well ordered. If $m$ is a
singular point, then necessarily $n_1$ is the \emph{rank} of $dF(m)$.
Then one can show the following lemma:
\begin{lemm}
  \label{lemm:equiv-forte}
  If $\mathbf{f}\seq\mathbf{g}$ at a point $m$, then there exist well
  ordered basis of $\mathbf{f}$ and $\mathbf{g}$ associated to a
  common partition $n=n_1+\cdots+n_d$, and an $n\times n$ matrix $N$
  with coefficients in $\C_\mathbf{f}$ ($=\C_\mathbf{g}$) such that
  $N(m)$ is block-diagonal (the i-eth block being of size $n_i$) and
  \[
  (g_1-g_1(m),\dots,g_n-g_n(m))=N\cdot(f_1-f_1(m),\dots,f_n-f_n(m)).
  \]
\end{lemm}

\subsection{Quantum mechanics}
The quantum analysis will be performed using $\h$-pseudo-differential
operators. This theory is now well established and very robust. It is
simultaneously a quantum theory dealing with self-adjoint operators on
a Hilbert space and a semiclassical theory, dealing with
$\h$-deformations of classical Hamiltonians. The reader is invited to
check the references \cite{robert}, \cite{dimassi-sjostrand} or
\cite{courscolin}.

Let $X$ be an $n$-dimensional smooth manifold equipped with a density
$|dx|$. The Hilbert space is the corresponding $L^2(X,|dx|)$.
$\h$-\pdos{} on $X$ act on $L^2(X,|dx|)$ and have symbols defined on
the cotangent bundle $M=T^*X$. Locally in $x$, any such operator can
be obtained from a \emph{full symbol} $a(x,\xi;\h)$ by the Weyl
quantisation formula
\begin{equation}
  (A u) (x) = (\textup{Op}^w_h(a) u) (x) =
  \frac{1}{(2\pi\h)^n}\int_{\RM^{2n}} e^{\frac{i}{\h}(x-y).\xi}
  a({\textstyle\frac{x+y}{2}},\xi;h) u(y)|dyd\xi|
  \label{equ:weyl}
\end{equation}
Such a full symbol is not invariant under coordinate changes, but if
we restrict to volume preserving diffeomorphisms, and to symbols
admitting an asymptotic expansion\footnote{Such semiclassical symbols
  are usually called ``classical''...} of the form
$\h^ka_0(x,\xi)+\h^{k+1}a_1(x,\xi)+O(\h^{k+2})$, then $k$ is called
the \emph{order} of the operator, and the \emph{principal symbol}
$a_0$ and the \emph{subprincipal symbol} $a_1$ are intrinsically
defined as functions on $T^*X$.

The space of $\h$-\pdos{} is a graded algebra.  We have a symbolic
calculus, which means that the product of operators of order zero is
still a \pdo{} of order zero whose principal symbol is the product of
the original principal symbols.  Moreover the commutation bracket of
operators of order zero is a \pdo{} of order 1 whose principal symbol
is $\frac{1}{i}$ times the Poisson bracket of the original principal
symbols.

By definition, a \emph{semiclassical completely integrable system} is
the data of $n$ self-adjoint $\h$-\pdos{} $P_1,\dots,P_n$ of order
zero which pairwise commute modulo $\oh$: $[P_i,P_j]=\oh$ and whose
principal symbols $p_j$ have almost everywhere independent
differentials. Then of course these symbols $(p_1,\dots,p_n)$ define a
classical completely integrable system on $T^*X$.

We are mainly interested in the microlocal behaviour of the $P_j$'s to
construct \emph{joint quasimodes}, ie. solutions $u$ of the system of
equations
\begin{equation}
  \label{equ:systeme}
  P_j u =\oh \qquad \forall j=1,\dots,n,
\end{equation}
where $u$ is a distribution on $X$ depending on $\h$ in a temperate
way (which we don't explicit here).  More generally we shall deal with
\emph{microlocal solutions} of this system. Roughly speaking $u$ is a
microlocal solution at a point $m\in T^*X$ if~\eqref{equ:systeme}
holds when both sides are multiplied on the left by a \pdo{} whose
principal symbol does not vanish at $m$.  The set of points in $T^*X$
where a distribution $u$ does not vanish microlocally is called the
\emph{microsupport} of $u$. A microlocal solution
of~\eqref{equ:systeme} has therefore a microsupport in the level set
$p^{-1}(0)$.

As quantum operators, the $P_j$'s have a spectrum, which we shall
always assume to be discrete (this is the case for instance when the
momentum map $p=(p_1,\dots,p_n)$ is proper). If they commute exactly:
$[P_i,P_j]=0$ then they also have a \emph{joint spectrum}, which is
the set of n-uples of eigenvalues
$(\lambda_1,\dots,\lambda_n)\in\RM^n$ associated to a common
eigenfunction. Then the construction of joint quasimodes goes a long
way in describing this joint spectrum modulo $\oh$.

In analogy with the classical case, for each semiclassical completely
integrable system $P_1,\dots,P_n$ and for any open subset $U\subset
M=T^*X$, we define the semiclassical momentum algebra
$\mathbf{P}_{\restr U}$ to be the linear span of the $P_j$'s, and the
semiclassical commutant $\C_\mathbf{P}(U)$ to be the Lie algebra of
all $\h$-\pdos{} commuting with $\mathbf{P}$ microlocally in $U$. As
before, $\mathbf{P}$ and $\C_\mathbf{P}$ are presheaves over $M$.

The weak and strong equivalences for semiclassical systems are defined
as in the classical case. For instance, we will use the following:
\begin{defi}
  Two semiclassical completely integrable systems with momentum
  algebra $\mathbf{P}$ and $\mathbf{Q}$ are called \textbf{strongly
    equivalent} at a point $m\in M$ if, microlocally near $m$,
\[
\C_{\mathbf{P}}.\left(\mathbf{P}-\mathbf{p}(m)\right) =
\C_{\mathbf{Q}}.\left(\mathbf{Q}-\mathbf{q}(m)\right).
\]
This will be denoted by $\mathbf{P}\seq\mathbf{Q}$.
\end{defi}
The ``quantised'' version of lemma~\ref{lemm:equiv-forte} then holds,
where $N$ becomes a matrix with pseudodifferential coefficients.
Notice that if $\mathbf{P}\seq\mathbf{Q}$ at $m$ then
$\mathbf{P}-\mathbf{p}(m)$ and $\mathbf{Q}-\mathbf{q}(m)$ share the
same microlocal joint quasimodes. In a sense, the set (presheaf) of
all joint quasimodes for a semiclassical system is the quantum
analogue of the classical Lagrangian foliation.

Given a semiclassical system $\mathbf{P}$, we have an underlying
classical system given by the principal symbols. The
\emph{subprincipal symbols} $r_j$ can then be seen as a small
deformation of the induced Lagrangian foliation; more precisely, they
define a family of 1-forms $\kappa_c$ on the leaves $p^{-1}(c)$ by the
formula
\begin{equation}
  \kappa_c(\ham{p_j})=-r_j,
  \label{equ:sous-principale}
\end{equation}
It is easy to check (using for instance the Darboux-Carathéodory
theorem below) that near any regular point of $p^{-1}(c)$, $\kappa_c$
is a smooth closed 1-form. It is called the \emph{subprincipal form}
of the semiclassical system.

\subsection{Canonical transformations}
The main strength of \pdos{} is the possibility of transforming them
according to any local symplectomorphism. If $\chi$ is a local
symplectomorphism of $T^*X$ near $m$ then there exist a \emph{\fio}
$U$ which is a bounded operator on $L^2(X)$ such that for any \pdo{}
$P$ with principal symbol $p$, the operator $U^{-1}PU$ is a
pseudodifferential operator whose principal symbol near $m$ is
$p\circ\chi^{-1}$ (Egorov theorem~\cite{robert,courscolin}). 

Constructing more global \fio{}s using a partition of unity is not
difficult, provided the following obstruction vanishes: let $\alpha$
be the canonical Liouville 1-form of $T^*X$. Then
$\alpha-\chi^*\alpha$ is closed. The obstruction is its cohomology
class. In other words $\chi$ should be ``exact'' in the sense that it
preserves integrals of $\alpha$ along closed loops.

If $\mathbf{P}$ is a quantum completely integrable system and $U$ a
\fio{} associated to a canonical transformation $\chi$ then
$U^{-1}\mathbf{P}U$ is a quantum completely integrable system with
momentum algebra $\mathbf{p}\circ\chi^{-1}$. Moreover the subprincipal
form $\kappa_c$ is modified only by the addition of an exact 1-form.

\paragraph{Disclaimer ---} I have deliberately included no example in
this review, for several reasons. One is brevity. Another one is that
the interested reader can find many examples in the bibliography. But
maybe most importantly one of the points of using theoretical normal
forms is to simplify the study; and it turns out that even for the
simplest examples the normal forms give a much easier way to discover
interesting phenomena than explicit calculations (which are
furthermore very often impossible). Nonetheless I am still convinced
that examples are essential, non only to motivate the theory, but also
to discover the features that, finally, may turn out to be easier to
cope with using the general theory...

\section{Local study}
The local behaviour of a completely integrable system can be very rich
and is far from being thoroughly understood in general. We review here
the current state of the art for the $\Cinf$ category and show how it
applies to quantum systems.

\subsection{Regular points}
Let $(f_1,\dots,f_n)$ be a classical completely integrable system on a
$2n$-symplectic manifold $M$, with momentum map $F$.
\begin{defi}
  A point $m\in M$ is called \textbf{regular} for $F$ if $dF(m)$ has
  maximal rank ($n$). In other words $df_1\wedge\cdots\wedge
  df_n(m)\neq 0$.
\end{defi}
By the local submersion theorem, the fibres $F^{-1}(c)$ for $c$ close
to $F(m)$ are locally $n$-dimensional submanifolds near a regular
point $m$. The local structure of regular points of completely
integrable systems is very simple. It is actually entirely described
by the following classical theorem:
\begin{theo}[Darboux-Carathéodory]
  If $m$ is regular, $F$ is symplectically conjugate near $m$ to the
  linear fibration $(\xi_1,\dots,\xi_n)$ on the linear symplectic
  space $\RM^{2n}$ with coordinates
  $(x_1,\dots,x_n,\xi_1,\dots,\xi_n)$ and symplectic form $\sum_i
  d\xi_i\wedge dx_i$.
  
  In other words, there exists functions $\phi_1,\dots,\phi_n$ on $M$
  such that
 \[
 (\phi_1, \dots, \phi_n, f_1, \dots, f_n)
  \] 
  is a system of canonical coordinates in a neighbourhood of $m$.
\end{theo}

In principle the name of Liouville should be associated with this
theorem, since well before Darboux and Carath\'eodory, Liouville gave
a very nice and explicit formula for the functions $\phi_j$. This
result published in 1855~\cite{liouville-note} explains the local
integration of the flow of any completely integrable Hamiltonian
(possibly depending on time) near a regular point of the foliation in
terms of the famous Liouville 1-form $\sum_i\xi_idx_i$. With this
respect it implies the Darboux-Carathéodory theorem, even if
Liouville' s formulation is more complicated.

The Darboux-Carathéodory theorem has some simple but very important
corollaries.
\begin{lemm}
  \label{lemm:regulier}
  Locally near a regular point, the commutant $\C_\mathbf{f}$ is the
  set of functions of the form $\phy(f_1,\dots,f_n)$, where
  $\phy\in\Cinf(\RM^n)$.
\end{lemm}
\begin{demo}
  Apply Darboux-Carathéodory.
\end{demo}
\begin{prop}
  \label{prop:commutative}
  $\C_\mathbf{f}$ is a commutative Lie-Poisson algebra.
\end{prop}
\begin{demo}
  As we saw already, it is a Lie algebra due to the Jacobi identity.
  It is a Poisson algebra due to the Leibniz identity. By the
  preceding lemma, it is commutative near regular points, and hence
  everywhere.
\end{demo}
As a consequence of this proposition, we have the following useful
lemma:
\begin{lemm}
  If $\mathbf{f}$ and $\mathbf{g}$ are two momentum algebras, then
  \[
  \mathbf{f}\sim\mathbf{g} \ssi \mathbf{f}\subset \C_\mathbf{g} .
  \]
\end{lemm}
\begin{coro}
  \label{coro:restriction}
  If $\mathbf{f}\sim_U\mathbf{g}$ and $V\subset U$ then
  $\mathbf{f}\sim_V\mathbf{g}$.
\end{coro}
Finally, lemma~\ref{lemm:regulier} implies the following
characterisation (which will not hold in the singular case)
\begin{prop}
  \label{prop:equiv-reg}
  If $m$ is regular for both $\mathbf{f}$ and $\mathbf{g}$ then
  $\mathbf{f}\sim\mathbf{g}$ near $m$ if and only if
  $(f_1,\dots,f_n)=\phy(g_1,\dots,g_n)$, where $\phy$ is a local
  diffeomorphism of $\RM^n$.
\end{prop}
In this case, the weak equivalence is identical to the strong one
(recall lemma~\ref{lemm:equiv-forte} above).

\paragraph{Semiclassics ---}
The Darboux-Carathéodory theorem admits a semiclassical analogue which
is again simple but powerful. The first proof in the framework of
homogeneous \pdos{} is due to Colin de Verdi\`ere~\cite{colinI}
(although the case $n=1$ was already treated by Duistermaat and
H\"ormander~\cite{FIO2}).

Let $P_1,\dots, P_n$ be a semiclassical completely integrable system
on $M=T^*X$ and let $p=(p_1,\dots,p_n)$ be the classical momentum map
consisting of the principal symbols. We also use $P=(P_1,\dots,P_n)$
for the quantum momentum map, or ``quantum fibration''.
\begin{theo}
  \label{theo:darboux-caratheodory-semicla}
  If $m$ is regular then $P$ is microlocally conjugate near $m$ to the
  fibration
  $(\frac{\h}{i}\deriv{}{x_1},\dots,\frac{\h}{i}\deriv{}{x_n})$ acting
  on $\RM^n$.
  
  In other words, there exists a \fio{} $U$ defined near $m$ and
  microlocally unitary such that
  $U^{-1}P_jU=\frac{\h}{i}\deriv{}{x_j}$.
\end{theo}
\begin{demo}
  Consider the symplectomorphism given by Darboux-Carathéodory's
  theorem and $U_0$ a \fio{} quantising it near $m$. The result is
  thus obtained modulo \pdos{} of order 1. To correct this error one
  conjugates again by a \pdo{} of the form $\exp(iA)$, where $A$ is a
  \pdo. Since we are only after a microlocal result (and hence modulo
  $\oh$), it is enough to show that there is a neighbourhood of $m$ in
  which one can solve to any order in $\h$, which is a simple
  exercise.
\end{demo}
\begin{rema}
  It is not necessary for $P$ to be self-adjoint. As long as the
  principal symbol is real valued, the theorem still holds, but the
  unitariness of $U$ is lost.
\end{rema}
Using this theorem one can check that all the corollaries of the
classical Darboux-Carathéodory theorem we have mentioned above still
hold in the semiclassical framework. It is even more important to see
that one can now describe microlocal joint quasimodes near regular
points.
\begin{prop}[\cite{san-focus}]
  \label{prop:solu_reg}
  If $m$ is a regular point, the space of microlocal solutions of the
  system
  \[
  P_j u = \ohb \quad \textrm{near } m \qquad \forall j=1,\dots,n
  \]
  is a $\CM_\h$-module of rank 1, generated by $U^{-1}\mathbf{1}$,
  where $U$ is a \fio{} as in
  Theorem~\ref{theo:darboux-caratheodory-semicla} and $\mathbf{1}$ is
  a wave function microlocally equal to 1 near $0\in\RM^{2n}$.
\end{prop}
Here $\CM_\h$ is the natural ring that acts on microlocal solutions;
it is the set of complex numbers depending in a temperate way on $\h$
(see \cite{san-focus}).

\subsection{Singular points}

The singularity theory of integrable systems is certainly not
completely understood, even for classical systems. I present here my
personal perception of it.

The study of singularities of integrable systems is fundamental for
various reasons. On the one hand, because of the way an integrable
system is defined : $n$ functions on a manifold, it is expected (apart
from exceptional cases) that singularities will necessarily occur. On
the other hand these functions define a dynamical system such that
their singularities correspond to fixed points and relative equilibria
of the system, which are of course one of the main characteristics of
the dynamics. From a semiclassical viewpoint, we know furthermore that
important wave functions such as eigenfunctions of the system have a
microsupport which is invariant under the classical dynamics;
therefore, in a sense that I shall not present here (one should talk
about semiclassical measures), they concentrate near hyperbolic
singularities (see for instance~\cite{colin-p} and the work of
Toth~\cite{toth-expected})). This concentration entails not only the
growth in norm of eigenfunctions (see for
instance~\cite{toth-zelditch-1}) but also a higher local density of
eigenvalues (see figure~\ref{fig:spec-hyp} in
section~\ref{sec:sing-fibres} below and the
articles~\cite{colin-p2,san-focus,san-colin}).

The singularities of a Hamiltonian system can be approached either
through the study of the flow of the vector fields --- this is the
``dynamical systems'' viewpoint --- or through the study of the
Hamiltonian functions themselves --- this is the ``foliation''
perspective. In the case of completely integrable systems, both
aspects are equivalent because the vector fields of the $n$ functions
$f_1,\dots,f_n$ form a basis of the tangent spaces of the leaves of
the foliation $f_i=\textup{const}_i$, at least for regular points. I
shall always tend to be on the foliation side, which displays more
clearly the geometry of the problem.

However, the foliations we are interested in are \emph{singular}, and
the notion of a singular foliation is already delicate. Generally
speaking these foliations are of Stefan-S\"u{\ss}mann
type~\cite{stefan}~: the leaves are defined by an integrable
distribution of vector fields. But they are more than that: they are
Hamiltonian, and they are \emph{almost regular} in the sense that the
singular leaves cannot fill up a domain of positive measure. The
precise way of dealing with these singular foliations is encoded in
the way we define two equivalent foliations; for us, this will be the
strong equivalence of definition~\ref{defi:eq-forte}.

\paragraph{Non-degenerate singularities ---}
In singularity theory for differentiable functions, ``generic''
singularities are Morse singularities. In the theory of completely
integrable systems there exists a very natural analogue of the notion
of Morse singularities (or more generally of Morse-Bott singularities
if one allows critical submanifolds). These so-called
\emph{non-degenerate} singularities are now well defined (and
``exemplified'') in the literature, so I will only recall briefly the
definition.

Let $F=(f_1,\dots,f_n)$ be a completely integrable system on $M$ with
momentum algebra $\mathbf{f}$.
\begin{defi}[\cite{vey}] A \emph{fixed} point $m\in M$ is
  called \textbf{non-degenerate} if the Hessians $d^2f_j(m)$ span a
  Cartan subalgebra of the Lie algebra of quadratic forms on $T_m M$
  (equipped with the linearised Poisson bracket).
\end{defi}
A friendlier characterisation is that a generic linear combination of
the linearised vector fields at $m$ (these are Hamiltonian matrices :
in $\textup{sp}(2n,\RM)$) should admit $2n$ distinct eigenvalues. The
definition above applies to a fixed point. But more generally if
$dF(m)$ has corank $r$ one can assume that $df_1(m),\dots,df_{n-r}(m)$
are linearly independent; then we consider the restriction of
$f_{n-r+1},\dots,f_n$ to the symplectic manifold $\Sigma$ obtained by
local symplectic reduction under the action of $f_1,\dots,f_{n-r}$. We
shall say that $m$ is non-degenerate (or \emph{transversally
  non-degenerate}) whenever $m$ is a non-degenerate fixed point for
this restriction of the system to $\Sigma$.

The linear approximation or \emph{linear model}\footnote{The term
  ``linear'' refers to the linearisation of the vector fields
  $\ham{f_i}$ at a fixed point; of course the functions $f_j$
  themselves do not become linear, but quadratic.} of such a critical
point is the system
$\mathbf{f_0}:=(\xi_1,\dots,\xi_{n-r},q_1,\dots,q_r)$ on
$T^*\RM^{n-r}\times T_m\Sigma$, where the $q_j$'s form a basis of the
aforementioned Cartan subalgebra.
\begin{theo}[Eliasson's theorem \cite{eliasson,eliasson-these}]
  Non-degenerate critical points are linearisable: there exists a
  local symplectomorphism $\chi$ in the neighbourhood of $m$ such that
\[
\chi^*\mathbf{f} \seq \mathbf{f_0} .
\] 
\end{theo}
In order to use this theorem one has to understand the linear
classification of Cartan subalgebras of $\textup{sp}(2n,\RM)$. This
follows from the work of Williamson~\cite{williamson}, which shows
that any such Cartan subalgebra has a basis build with three type of
blocks: two uni-dimensional ones (the elliptic block: $q=x^2+\xi^2$
and the real hyperbolic one: $q=x\xi$) and a two-dimensional block
called focus-focus or loxodromic or complex hyperbolic:
$q_1=x\xi+y\eta$, $q_2=x\eta-y\xi$. Notice that over $\CM$ the
classification is trivial since everything can be conjugate to the
``hyperbolic'' case $x\xi$. The analytic case of Eliasson's theorem
was proved by R\"u{\ss}mann~\cite{russmann} for two degrees of freedom
systems and by Vey~\cite{vey} in any dimension. In the $\Cinf$
category the \emph{lemme de Morse isochore} of Colin de Verdi\`ere and
Vey~\cite{colin-vey} implies Eliasson's result for one degree of
freedom systems. Eliasson's proof of the general case was somewhat
loose at a crucial step, but this has been recently completely
clarified~\cite{san-miranda}.

The strong equivalence relation for non-degenerate singularities is
equal to the weak equivalence and is fairly well understood. In
particular in case of real hyperbolic blocks it does not imply the
functional equivalence (which would be
``$\chi^*\mathbf{f}=\phy\circ\mathbf{f_0}$'')) while this is indeed
the case otherwise. For more details refer to~\cite{san-fn}.

\paragraph{Semiclassics ---}
It is possible to prove a semiclassical version of Eliasson's
theorem, which is crucial for further study of singular
Bohr-Sommerfeld conditions as in~\cite{san-focus} or for the estimate
by Zelditch and Toth of norms of
eigenfunctions~\cite{toth-zelditch-2}. The new semiclassical feature
is the appearance of formal series in $\h$ of microlocal invariants.
\begin{theo}[\cite{san-fn}]
  \label{theo:fn-semicla}
  If $m$ is non-degenerate of corank $r$, there exists a \fio{} $U$
  defined near $m$ and microlocally unitary and there exists formal
  series $\alpha_j(\h)\in\CM\formel{\h}$, $j=1,\dots,r$ such that
\[
  U \left(\begin{array}{c} P_1 \\ \vdots \\ P_n
    \end{array}\right)U^{-1}
  \seq \left(
    \begin{array}{c}
      \hat{\xi}_1 \\ \vdots \\ \hat{\xi}_{n-r} \\
      \hat{q}_1-\h\alpha_1(\h) \\ \vdots \\ \hat{q}_r-\h\alpha_r(\h)
    \end{array}
  \right)
  \]
  acting on $\RM^{n-r}\times\RM^r$, microlocally near $m$.
\end{theo}
In this statement we have used the hat for standard Weyl quantisation
(for instance $\hat{\xi}_j=\frac{\h}{i}\deriv{}{x_i}$ and
$\widehat{x\xi}=\frac{\h}{i}(x\deriv{}{x} + \frac{1}{2})$). The
introduction of the series $\alpha_j(\h)$ is necessary to go from the
weak equivalence to the strong one.

This theorem is well adapted to the microlocal resolution of the
system $P_ju=\oh$ since the latter is transformed into the system
$(\hat{q}_j-\alpha_j)u=\oh$, which can be solved explicitly. Since
the model system is uncoupled, one just has to study separately each
block, and one can show the following facts: for an elliptic block,
the space of microlocal solutions (in the sense of
proposition~\ref{prop:solu_reg}) has dimension 1; for a real
hyperbolic block, it has dimension 2~\cite{colin-p2}; for a
focus-focus block, it has dimension 1~\cite{san-focus}.

\subsection{More degenerate singularities}

For the moment very little is known concerning degenerate
singularities. The most natural approach seems to be via algebraic
geometry, as in \cite{garay,garay-straten}. For one degree of freedom
analytic systems, a more concrete method is presented in
\cite{colin-singularites}, which explicitly displays the relevant
versal unfoldings. For a general linearisation result in the analytic
category, see also \cite{zung-birkhoff}. I am not aware of similar
results in the $\Cinf$ category.

\section{Semi-global study}

If one aims at understanding the classical geometry of a completely
integrable foliation or its microlocal analysis, the
\emph{semi-global} aspect is probably the most fundamental. The
terminology semi-global refers to anything that deals with an
invariant neighbourhood of a leaf of the foliation. This semi-global
study is what allows for instance the construction of quasimodes
associated to a Lagrangian submanifold. Sometimes semi-global merely
reduces to local, when the leaf under consideration is a critical
point with only elliptic blocks.

\subsection{Regular fibres}

The analysis of neighbourhoods of regular fibres, based on the
\aangles{} theorem (also known as action-angle theorem) is now routine
and fully illustrated in the literature, for classical aspects as well
as for quantum ones. It is the foundation of the whole modern theory
of completely integrable systems (in the spirit of Duistermaat's
article~\cite{duistermaat}) but also of KAM-type perturbation
theorems. The microlocal analysis of action-angle variables starts
with the work of Colin de Verdi\`ere~\cite{colinII}, followed in the
$\h$ semiclassical theory by Charbonnel~\cite{charbonnel}, and more
recently by myself and various articles by Zelditch, Toth, Popov,
Sj\"ostrand and many others. The case of compact symplectic manifolds
has recently started, using the theory of Toeplitz
operators~\cite{charles-bs}.

Let $(f_1,\dots,f_n)$ be an integrable system on a symplectic manifold
$M$. In the rest of this article we shall always assume the momentum
map $F$ to be \emph{proper}: all fibres are compact. Let $c$ be a
regular value of $F$. If we restrict to an adequate invariant open
set, we can always assume that the fibres of $F$ are connected. Let
$\Lambda_c:=F^{-1}(c)$. Fibres being compact and parallelisable (by
means of the vector fields $\ham{f_i}$), they are tori.
\begin{theo}[\aangles]
  If $\Lambda_c$ is regular, there exists a symplectomorphism $\chi$
  from $T^*\T^n$ into $M$ sending the zero section onto $\Lambda_c$
  such that 
 \[
  \chi^*\mathbf{f} \sim \mathbf{f_0},
  \]
  where $\mathbf{f_0}$ is the linear system $(\xi_1,\dots,\xi_n)$ on
  $T^*\T^n$.
\end{theo}
Here and in what follows we identify $T^*\T^n$ with $\T^n\times\RM^n$
(where $\T=\RM/\ZM$) equipped with coordinates
$(x_1,\dots,x_n,\xi_1,\dots,\xi_n)$ such that the canonical Liouville
1-form is $\sum_i \xi_i dx_i$.  It is easy to see that the theorem
actually implies $\chi^*\mathbf{f}=\phy(\mathbf{f_0})$ for $\phy$ a
local diffeomorphism of $\RM^n$; this is usually the way the result is
stated. It is important to remark that $d\phy$ is an \emph{invariant}
of the system since it is determined by \emph{periods} of periodic
trajectories if the initial system.  Regarded as functions on $M$ the
$\xi_j$'s are called \emph{actions} of the system for one can find a
primitive $\alpha$ of $\omega$ in a neighbourhood of $\Lambda_c$ such
that the $\xi_j$'s are integrals of $\alpha$ on a basis of cycles of
$\Lambda_c$ depending smoothly on $c$.

For semiclassical purposes, one needs an ``exact'' version of the
\aangles{} theorem in the sense that, given a primitive $\alpha$ of
$\omega$, a symplectomorphism is called exact when it preserves the
integrals of $\alpha$ along closed paths (action integrals). We get
immediately:
\begin{theo}
  \label{theo:action-angle-exact}
  If $\Lambda_c$ is regular, there exists an exact symplectomorphism
  $\chi$ from $T^*\T^n$ into $M$ sending the section
  $(\xi_1,\dots,\xi_n)=(a_1,\dots,a_n)=\textrm{const}$ onto
  $\Lambda_c$ such that
  \[
  \chi^*\mathbf{f} \sim \mathbf{f_0}
  \]
  if and only if $a_i=\int_{\gamma_i}\alpha$, where $\gamma_i$ is the
  cycle on $\Lambda_c$ corresponding via $\chi$ to the i-eth canonical
  cycle on $\T^n$.
\end{theo}

\paragraph{Semiclassics ---} Let $(P_1,\dots,P_n)$ be a semiclassical
completely integrable system whose principal symbols define a proper
momentum map. We still denote by $\Lambda_c$ the Lagrangian leaves.
$\alpha$ is the canonical 1-form of $M=T^*X$ and as before we let
$a=(a_1,\dots,a_n)$ be the action integrals along the basis of cycles
of $\Lambda_c$ defined by the chosen actions $\xi_j$.
\begin{theo}[Semiclassical action-angle]
  If $\Lambda_c$ is regular there exists formal series
  $\lambda_j(\h)\in\CM\formel{\h}$ and a \fio{} $U$ associated to an
  exact symplectomorphism from $T^*\T^n$ to $M$ sending the
  ``$\xi=a$'' section onto $\Lambda_c$ such that, microlocally near
  the ``$\xi=a$'' section, we have
  \[
  U(P_1,\dots,P_n)U^{-1} \seq
  (\hat{\xi_1}-\h\lambda_1(\h),\dots,\hat{\xi_n}-\h\lambda_n(\h)),
  \]
  acting on $\T^n$.
\end{theo}
Explicitly this means
\[
U(P_1-p_1(m),\dots,P_n-p_n(m)) U^{-1} =
N\cdot(\hat{\xi_1}-\tilde{\lambda}_1(\h),\dots,
\hat{\xi_n}-\tilde{\lambda}_n(\h)),
\]
where $m$ is any point of $\Lambda_c$, $N$ is an $n\times n$
microlocally invertible matrix of \pdos{} and
$\tilde{\lambda}_j(\h)=a_j+\h\lambda_j(\h)$. One can see that the
action $a_j$ can be considered as the first semiclassical invariant.
The second term is given by integrals of the \emph{subprincipal form}
(see definition~\eqref{equ:sous-principale}
page~\pageref{equ:sous-principale}) and of the \emph{Maslov cocycle} on
the Lagrangian $\Lambda_c$ (see also~\cite{san-focus}).

The following statement is a direct consequence of the theorem.
\begin{theo}[Regular Bohr-Sommerfeld quasimodes~\cite{san-focus}]
  \label{theo:quasi-modes}
  There is a non trivial microlocal solution of the system $P_ju=\ohb$
  (which is therefore microlocalised on $\Lambda_0$) if and only if
  $\tilde{\lambda}_j(\h)\in 2\pi\h\ZM$. The solution is unique (in the
  sense of proposition~\ref{prop:solu_reg}).
\end{theo}
From this one can deduce the \emph{regular Bohr-Sommerfeld conditions}
in the following way. Suppose $\mathbf{P}=\mathbf{P}^E$ depends on a
parameter $E\in\RM^n$ in such a way that for any $m$ near
$\Lambda_0^E=(p^E)^{-1}(0)$, the principal symbols map
$(E_1,\dots,E_n)\fleche (p_1^E,\dots,p_n^E)(m)$ is a local
diffeomorphism.
\begin{defi}
  We call \textbf{microlocal joint spectrum} the set
  $\Sigma_{\h}(P_1^E,\dots,P_n^E)$ of all $E\in\RM^n$ such that the
  system $P_j^Eu_{\h}=\ohb$, $j=1,\dots,n$ admits a non trivial
  microlocal solution on the whole fibre $\Lambda_0^E$.
\end{defi}
The typical case if of course $P_j^E=P_j-E_j$. From our perspective it
is often wiser to forget the linear dependence on $E$, which is not
invariant under strong equivalence.

The regular Bohr-Sommerfeld conditions are obtained from
theorem~\ref{theo:quasi-modes} if one remarks that it still holds
``with parameters''. They can be stated as follows:
\begin{theo}[\cite{san-focus}]
  \label{theo:BS-reg}
  The microlocal joint spectrum consists of solutions $E$ of
  \[
  \tilde{\lambda}^E_j(\h)\in 2\pi\h\ZM,
  \]
  where
  \begin{equation}
    \tilde{\lambda}^E_j(\h) = \int_{\gamma_i^E}\alpha +
    \h\int_{\gamma_i^E}\kappa^E +\h\frac{\mu(\gamma^E)\pi}{2}
    + O(\h^2).
    \label{equ:BS-reg}
  \end{equation}
  Here $\kappa^E$ is the subprincipal 1-form on $\Lambda_0^E$ and
  $\mu$ the Maslov cocycle. $(\gamma_1^E,\dots\gamma_n^E)$ is any
  basis of cycles on $\Lambda_0^E$.
\end{theo}

\subsection{Singular fibres}
\label{sec:sing-fibres}

This section is devoted to the semi-global structure of fibres with
non-degenerate singularities. I am not aware of any semi-global result
for more degenerate singularities. The topological analysis of
non-degenerate singular fibres was mainly initiated by
Fomenko~\cite{fomenko}, and successfully expanded by a number of his
students. See~\cite{bolsinov-fomenko-book}.

\paragraph{Elliptic case ---} 
Near an elliptic fixed point, the fibres are small tori and are
entirely described by the local normal form, for classical systems as
well as for semiclassical ones (the system is reduced to a set of
uncoupled harmonic oscillators). Therefore I shall not talk about this
type of singularity any further... even if strictly speaking the
semi-global semiclassical study has not been fully carried out for
\emph{transversally elliptic} singularities. But no particular
difficulties are expected in that case.

\paragraph{Focus-focus case ---}
Eliasson's theorem gives the local structure of focus-focus
singularities. Several people have noticed (in the years 1996-1997)
that this was enough to determine the \emph{monodromy} of the
foliation around the singular fibre; I'll expand on this in
section~\ref{sec:global}.  Actually this local structure is a starting
point for understanding much more: the semi-global classification of a
singular fibre of focus-focus type. Unlike monodromy which is a
topological invariant, already observed in torus fibrations without
Hamiltonian structure, the semi-global classification involves purely
symplectic invariants.

Let $F=(f_1,f_2)$ be a completely integrable system with two degrees
of freedom on a 4-dimensional symplectic manifold $M$. Let $m$ be a
critical point of focus-focus type; we assume for simplicity that
$F(m)=0$, and that the (compact, connected) fibre $\Lambda_0$ does not
contain other critical points. One can show that $\Lambda_0$ is a
``pinched'' torus (Lagrangian immersion of a sphere $S^2$ with a
transversal double point), surrounded by regular fibres which are
standard $\T^2$ tori. What are the semi-global invariants associated
to this singular fibration ?

\begin{figure}[htbp]
  \begin{center}
    \includegraphics{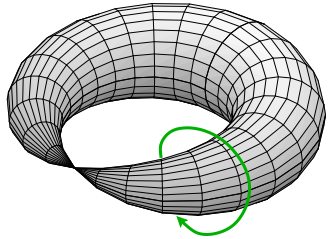}
    \caption{Vanishing cycle for the pinched torus}
    \label{fig:evanescent}
  \end{center}
\end{figure}

One of the major characteristics of focus-focus singularities is the
existence of a \emph{Hamiltonian action of $S^1$} that commutes with
the flow of the system, in a neighbourhood of $\Lambda_0$. Indeed, let
us start by applying Eliasson's theorem near $m$ to reduce to a
momentum map $F=(f_1,f_2)$ which is equal near $m$ to the canonical
focus-focus basis $(x\xi+y\eta,x\eta-y\xi)$. Then $f_2$ is the
periodic Hamiltonian we are looking for; it can also be identified
with an action integral associated to the vanishing cycle of the
pinched torus (cf.  fig.~\ref{fig:evanescent}).

Let $c$ be a regular value for $F$, close to 0. Given a point $A$ on
$\Lambda_c$, we define $\tau_1(c)>0$ to be the time of first return
for the $\ham{f_1}$-flow on the orbit of $A$ under the
$\ham{f_2}$-flow, and $\tau_2(c)\in\RM/2\pi\ZM$ to be the time it
takes to return to $A$ under the action of $\ham{f_2}$. Let us define
\[
\sigma_1(c)=\tau_1(c)+\re(\ln c) \quad \textrm{ et } \quad
\sigma_2(c)=\tau_2(c)-\im(\ln c),
\]
where we identified $c=c_1+ic_2$. One shows that
$\sigma:=\sigma_1(c)dc_1+\sigma_2(c)dc_2$ is a closed $\Cinf$ 1-form
in a neighbourhood of this origin. Let $s$ be the primitive of $\sigma$
vanishing at the origin. Let $[\mathbf{f]}$ be the foliation
associated to the system in a neighbourhood of $\Lambda_0$ (ie. the
equivalence class of $\mathbf{f}$ modulo strong equivalence). Let
$S([\mathbf{f]})$ be the Taylor expansion of $s$.
\begin{theo}[\cite{san-semi-global}]
\label{theo:invariants}
$S([\mathbf{f]})$ completely characterises the singular foliation in a
neighbourhood of $\Lambda_0$, which means:
\begin{itemize}
\item $S([\mathbf{f]})$ is well-defined (it does not depend on the
  choice of Eliasson's local chart and is invariant under strong
  equivalence);
\item If $[\mathbf{f]}$ and $[\mathbf{g]}$ are two singular foliations
  in a neighbourhood of focus-focus leaves, and satisfy
  $S([\mathbf{f}])=S([\mathbf{g}])$, then there exists a semi-global
  symplectomorphism $\chi$ such that $\chi^*\mathbf{f}\seq\mathbf{g}$;
\item If $T$ is any formal series in $\RM\formel{X,Y}$ without
  constant term, then there exists a singular foliation of focus-focus
  type $\mathbf{f}$ such that $T=S([\mathbf{f}])$.
\end{itemize}
\end{theo}
\begin{rema}
  The fact that two focus-focus fibrations are always semi-globally
  \emph{topologically} conjugate was already proved by
  Zung~\cite{zung-I}, who introduced various topological notions of
  equivalence. In our language this means that the topological class
  is invariant by strong equivalence. The theorem shows that this is
  no longer the case for the symplectic class.
\end{rema}
\begin{rema}
  \label{rema:regularisation}
  $S$ can be interpreted as a \emph{regularised} (or desingularised)
  \emph{action}. Indeed if $\gamma_c$ is the loop on $\Lambda_c$
  defined just as in the description of $\tau_j$ above, and if
  $\alpha$ is a semi-global primitive of the symplectic form $\omega$,
  let $\mathcal{A}(c)=\int_{\gamma_c}\alpha$; then
  \parbox{0.9\textwidth}{ \[ S(c) = \mathcal{A}(c) - \mathcal{A}(0) +
    \re(c\ln c - c).
  \]}
\end{rema}

The semiclassical study of focus-focus fibres was carried out in the
article~\cite{san-focus}. \emph{Singular Bohr-Sommerfeld conditions}
are stated there which allow a complete description of the microlocal
joint spectrum in a neighbourhood of the critical value of the
momentum map. Unlike the case of standard action-angle variables, it
is not possible to proceed merely by applying a semiclassical
semi-global normal form; indeed, the classification theorem does not
provide us with an explicit model. And all the examples that I know
of, that can reasonably claim to be a ``typical model of focus-focus
singularity'', are not explicitly solvable. The strategy of that
paper is to see the microlocal solutions as global sections of a sheaf
which, because of the local normal forms, is a locally flat bundle;
then such a global section exists if and only if the \emph{holonomy}
of the sheaf is trivial. This approach may of course be used in the
regular case as well; the phase of this holonomy is then identified to
the \emph{semiclassical action integrals} of
formula~\eqref{equ:BS-reg}. In the singular case the adequate holonomy
is a \emph{regularisation} of the usual semiclassical action, in the
sense of remark~\ref{rema:regularisation} above. As a matter of fact
the first term of this holonomy is exactly the classical invariant
described in theorem~\ref{theo:invariants} above, which entails that
the symplectic equivalence class of the foliation is a \emph{spectral
  invariant} of the quantum system. See \cite{san-focus,san-X} or more
details.

\ouf

Besides semiclassics, theorem~\ref{theo:invariants} leads to a number
of various applications. One can for instance exploit the fact that
the set of symplectic equivalence classes of these foliations acquires
a vector space structure. That is what Symington does
in~\cite{symington-gen} to show that neighbourhoods of focus-focus
fibres are always symplectomorphic (after forgetting the foliation, of
course). For this one introduces functions $S_0$ and $S_1$ whose
Taylor expansions give the invariants of the two foliations, and
constructs a ``path of foliations'' by interpoling between $S_0$ and
$S_1$. Then a Moser type argument yields the result (since the
symplectic forms are cohomologous).

The theorem is also very useful for doing almost explicit calculations
in a neighbourhood of the fibre. For instance it is possible in this
way to determine the validity of non-degeneracy conditions that appear
in KAM type theorems\footnote{A nice discussion about theses various
  conditions can be found in~\cite{roy-conditions}}, for a
perturbation of a completely integrable system with a focus-focus
singularity (see also \cite{zung-kolmogorov}).
\begin{theo}[\cite{san-dullin}]
  Let $H$ be a completely integrable Hamiltonian with a loxodromic
  singularity at the origin (ie. $H$ admits a singular Lagrangian
  foliation of focus-focus type at the origin). Then in a
  neighbourhood of 0,
  \begin{itemize}
  \item Kolmogorov non-degeneracy condition if fulfilled on \emph{all
      tori} close to the critical fibre;
  \item the ``isoenergetic turning frequencies'' condition is
    fulfilled \emph{except} on a 1-parameter family of tori
    corresponding to a curve through the origin in the image of the
    momentum map which is transversal to the lines of constant energy
    $H$.
  \end{itemize}
\end{theo}

\begin{figure}[htbp]
  \begin{center}
    \input{huit.pstex_t}
    \caption{Leaf space of a real hyperbolic foliation}
    \label{fig:huit}
  \end{center}
\end{figure}
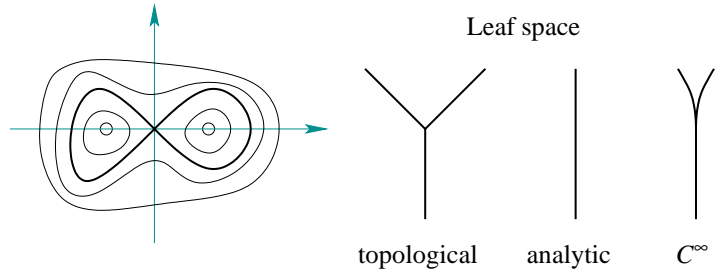

\paragraph{Hyperbolic case ---}
Just as elliptic blocks, hyperbolic blocks have dimension 1 (normal
form $q_i=x_i\xi_i$); however they turn out to be more complicated and
display a richer structure, due to two main reasons. The first one is
fundamental: the singular fibres are not localised near the
singularity; instead they consist of stable and unstable manifolds
that can connect several singular points. the second reason, more
technical, is that the natural $\Cinf$ structure of the space of
leaves is more involved. For instance in the case of the figure ''8''
(fig.~\ref{fig:huit}) the ``topological'' space of leaves is a ``Y'';
nevertheless from the real analytic viewpoint (when $H$ is analytic)
the space of leaves is just an interval (functions that commute with
$x\xi$ are functions of $x\xi$). In the $\Cinf$ category the space of
leaves is still a ``Y'', but whose hands have all derivatives equal at
the branching point; $\Cinf$ functions that commute with $x\xi$ are
locally described by two smooth functions $f^+(x\xi)$ and $f^-(x\xi)$
(for instance associated to the half-spaces $\pm x>0$) such that
$f^+-f^-$ is flat at the origin.

Concerning the semi-global aspect, the classification of hyperbolic
foliations has been carried out only for 1 degree of freedom, in
Toulet's thesis~\cite{toulet-these, dufour-mol-toul}. As a matter of
fact, it is possible to give a proof of the statement provided in that
note with similar methods to those that were used
in~\cite{san-semi-global}. In a slightly weaker context (orbital
equivalence) Bolsinov~\cite{bolsinov-orbital} has studied the case of
transversally hyperbolic singularities (codimension 1) in two degrees
of freedom.

For one degree of freedom systems, the critical fibre is a
\emph{graph} whose vertices have degree 4 (if all singularities of the
fibre are hyperbolic). The invariant is the graph itself, properly
\emph{decorated}. The corresponding semiclassical analysis was treated
by Colin de Verdi\`ere and Parisse~\cite{colin-p,colin-p2,colin-p3}.
The authors use the graph of Dufour-Molino-Toulet to read the correct
quantisation conditions.  Here again, a holonomy method is employed.

\begin{figure}[htbp]
  \begin{center}
    \includegraphics[angle=270, width=8cm]{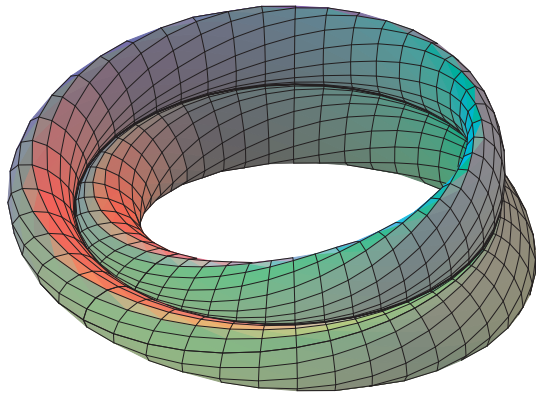}
    \caption{Hyperbolic critical fibre with  $\ZM/2\ZM$ symmetry (case
      of the $1:2$ resonance) }
    \label{fig:tore-double}
  \end{center}
\end{figure}

In the article~\cite{san-colin} we have studied in details the case of
two degrees of freedom with transversally hyperbolic singularities.
The set of critical points in the critical fibre is a union of
circles. In a neighbourhood of each circle, we reduce to a model
situation which may have a $\ZM/2\ZM$ symmetry (this is the case for
instance of the Birkhoff normal form in $1:2$ resonance, for a non
vanishing energy; the critical fibre is displayed in
figure~\ref{fig:tore-double}).
As for the unidimensional case, we construct a graph whose homology
serves to state the singular Bohr-Sommerfeld conditions. This graph is
abstractly the \emph{reduction} of the critical fibre by an $S^1$
action that we construct and which leave the foliation invariant. The
``subtlety'' is that critical circles with non-trivial $\ZM/2\ZM$
symmetry become vertices of \emph{degree 2} (instead of degree 4). On
the other hand the issue of the delicate $\Cinf$ structure on the
graph cannot be avoided either and leads to a somewhat involved proof
of the validity of these Bohr-Sommerfeld conditions. (This difficulty
was avoidable in dimension 1.)

\begin{figure}[htbp]
  \begin{center}
    \input{jspec-hyp.pstex_t}
    \caption{Part of the joint spectrum for a transversally hyperbolic
      singularity with $\ZM/2\ZM$ symmetry. Case of the $1:2$
      resonance. The horizontal axis carries the energy values and the
      vertical axis the values of the additional integral $K$.}
    \label{fig:spec-hyp}
  \end{center}
\end{figure}
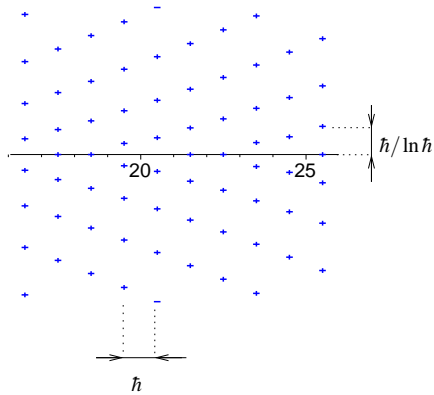

In this way one obtains in a very precise way the \emph{universal
  behaviour} of the microlocal joint spectrum near a transversally
hyperbolic separatrix, which yields amongst others the calculation of
the local density of eigenvalues (Weyl type formulas). As expected,
the distance between two points in the joint spectrum is of order $\h$
in the direction given by the periodic Hamiltonian, while it is of
order $\h/|\ln \h|$ in the transversal direction (cf.
fig.~\ref{fig:spec-hyp}).

By looking at a picture similar to figure~\ref{fig:spec-hyp},
Sadovski\'{\i} et Zhilinski\'{\i} had the idea that one could define
(and calculate) a \emph{fractional
  monodromy}~\cite{sadovski-zhilinski-nekh}, that reflects the
homology with rational coefficients of the singular torus fibration.
Using our work one should be able to write it in a rigorous way for a
general transversally hyperbolic singularity and calculate it by means
of the graph of the foliation. The ``singular rational affine
structure'' of the basis should be closely compared to the joint
spectrum as well.

\paragraph{Remaining cases ---}
To complete the study of two degrees of freedom integrable systems one
should still include the elliptic-hyperbolic case and the
hyperbolic-hyperbolic case.

The case of a system with a critical point splitting into an elliptic
block and a hyperbolic block is probably the more simple. It can be
regarded as a limit case of a transversally hyperbolic case whose
critical circles degenerate into a point, excluding the possibility of
a $\ZM/2\ZM$ symmetry. The critical fibre is therefore a graph of
degree 4 embedded in $M$.

Concerning semiclassics, one should obtain two quantisation
conditions: one related to the cohomology of the graph; the other to
the vanishing cycle.

\ouf

The hyperbolic-hyperbolic case is certainly a very good open
problem. It can be seen as a branching of four transversally
hyperbolic singularities, possibly coupled. The various possible
topologies are described in \cite{lerman}.

\paragraph{Higher degrees of freedom ---}

In principle no essential difficulty should arise in higher degrees of
freedom with non-degenerate singularities, since all the building
blocks have at most two degrees of freedom. I still believe this
should be very interesting to explore. On the classical side, the
topological theory was written by Zung~\cite{zung-I,zung-II}, but it
is not sufficient to turn it into semiclassics. On the other hand the
generalisation of the normal form near transversally hyperbolic orbits
(with discrete symmetry) has been recently achieved in
\cite{miranda-zung} and should be well suited for semiclassical
purposes.

\section{Global study}
\label{sec:global}

Up to now, we have considered various properties of integrable systems
in small neighbourhoods of minimal invariants objects: orbits of
points in $M$. How to get a global picture from them ? The adjective
``global'' covers several aspects, some qualitative and some
quantitative.

For instance, the aim of Duistermaat was to globalise properties given
by the \aangles{} theorem. He was thus interested in the fibration
over the set of \emph{regular points}, analysing the role of
monodromy, Chern class, and cohomology class of the symplectic form.

On the other hand the most natural way of globalising is to look for a
description of the symplectic manifold relying on formulas that
provide a ``localisation'' of global objects on \emph{singularities}
of the system --- just as Morse theory. Under the non-degeneracy
hypothesis for critical points, Zung made a thorough study and
displayed the importance of the \emph{integral affine structure} on
the base of the fibration~\cite{zung-II}.

This integral affine structure is another angle for dealing with the
global problem. In the most simple case of \emph{toric} completely
integrable systems (those whose flow defines an effective action of
$\T^n$) one can completely characterise the system consisting of the
symplectic manifold $M$ and the momentum map $F$ by means of the
image of $F$ which, in the integral affine manifold $\RM^n$, is a
convex rational polytope (Delzant's theorem~\cite{delzant}). It seems
that a natural way of generalising toric  systems is to allow only
non-degenerate singularities of elliptic or focus-focus types. These
systems are called \emph{almost toric} (the terminology was probably
introduced for the first time by Symington).

From the semiclassics viewpoint, the ``globalisation'' may refer to
the semiclassical analogues of the above geometrical globalisations.
It is also natural to consider the issue of passing from the
microlocal to the ``exact'': how to use microlocal constructions to
obtain results concerning the ``true'' Schrödinger operator acting on
the ``true'' Hilbert space $L^2(X)$ ? What is the relationship between
the microlocal spectrum and the exact spectrum ?

\subsection{The exact spectrum}

On $M=T^*X$ we are given a quantum completely integrable system
$P=(P_1,\dots,P_n)$.
\begin{defi}
  The \textbf{joint spectrum} of $P$ is the set of all
  $(E_1,\dots,E_n)\in\RM^n$ such that there exists a normalised
  element $\psi\in L^2(X)$ such that
 \[
  \forall i, \qquad P_i\Psi=E_i\Psi.
  \]
\end{defi}

The passage from microlocal to exact is based on two points. The
first one is very general and the second one is more adapted to our
vision of integrable systems.

The first point deals with the geometry and analysis at infinity: the
part of the Lagrangian foliation under study should be well separated
from further possible connected components of the fibration
$p=(p_1,\dots,p_n)$. This shall be granted by the assumption of local
properness of $p$: there exists a compact $K\subset\RM^n$ such that
$p^{-1}(K)$ is compact. If $X$ is compact, or $X=\RM^n$ and the $p_j$
have a good behaviour at infinity (for example all their derivatives
are bounded by a weight function in the sense of
H\"ormander~\cite{hormander-weyl}), one can show that in any compact
$K'$ whose interior is inside $K$, the joint spectrum is discrete (of
finite multiplicity)~\cite{charbonnel}.

The second point relies on the construction of microlocal quasimodes
and their microlocal \emph{multiplicity}, as in
proposition~\ref{prop:solu_reg}. The microlocal uniqueness of
solutions of the system not only shows that these quasimodes are good
approximations of eigenfunctions but also demonstrates that they form
a complete system, since they are microlocally orthogonal to each
other~\cite{san-focus}. This ensures that the microlocal spectrum is
really a perturbation of order $\oh$ of the exact spectrum, including
multiplicities.

\subsection{Regular fibrations: the case of monodromy}

I recall here the definition of monodromy and its semiclassical
consequences. I shall be very brief, referring for instance
to~\cite{san-mono,san-focus,san-dijon} for more details.

The \aangles{} theorem defines actions variables in a neighbourhood of
regular values of the momentum map $F$. Seen as local charts for the
open set $B_r$ of all regular values of $F$, they endow $B_r$ with the
structure of an integral affine manifold with structure group the
affine group $\textup{GL}(n,\ZM)\ltimes\RM^n$. By definition, the
\emph{affine monodromy} of the system is the holonomy of this affine
structure. Another way of defining an integral affine structure is to
specify a distribution of \emph{lattices} of maximal rank in each
tangent space. This lattice is the dual of the \emph{period lattice},
which is the set of all $(\tau_1,\dots,\tau_n)$ such that the
Hamiltonian vector field $\tau_i\ham{f_i}$ is $1$-periodic. The linear
part of the affine monodromy, simply called monodromy, is also the
holonomy of the flat bundle of homology groups of the fibres of the
torus fibration over $B_r$.

Following an idea of Cushman and Duistermaat, I have shown how this
monodromy can be read off from the microlocal \emph{joint spectrum} of
a corresponding quantum system. One just has to apply the regular
Bohr-Sommerfeld conditions which locally describe the joint spectrum
as a part of a lattice of type $\h\ZM^n$, whose mesh size tends to 0.
One can then define an asymptotic integral affine structure and show
that for $\h$ small enough it coincides with the classical affine
structure on $B_r$~\cite{san-mono}.

\subsection{Focus-focus and monodromy}

A remarkable properties of focus-focus singularities for a two degrees
of freedom system is that they imply the presence of a universal
nontrivial monodromy for the regular fibration around the critical
fibre. This results holds in general for topological torus fibrations
with a generic isolated critical fibre~\cite{moishezon,matsumoto}; in
the Hamiltonian case it was rediscovered by Nguyên Tiên
Zung~\cite{zung-focus}, Matveev~\cite{matveev}, and some others. The
particular feature of the Hamiltonian situation is that the monodromy
is oriented: while, in the topological case, the monodromy matrix is
$\left(
  \begin{array}[c]{cc}
1 & 0\\\pm 1 & 1
  \end{array}
\right)$, in the Hamiltonian case the sign is prescribed, and always
positive. Indeed if one chooses an orientation of $\RM^2$ it induces
through the momentum map $F$ and the symplectic form a natural
orientation on each Lagrangian torus, and hence on their
homology~\cite{cushman-san}. Concerning the quantum case,
Zhilinski\'{\i} suggests that the sign should be interpreted as the
fact that the ``lattice'' of eigenvalues around a focus-focus critical
value has a \emph{point defect} (in the sense that a number of points
was removed, and not added). This assertion can be verified by the
singular Bohr-Sommerfeld conditions which give the precise position of
eigenvalues near the singularity~\cite{san-focus}.
\begin{figure}[htbp]
  \begin{center}
    \input{monodromie.pstex_t}
    \caption{A joint spectrum with monodromy, and the calculation of
      the latter as an asymptotic affine holonomy}
    \label{fig:monodromie}
  \end{center}
\end{figure}
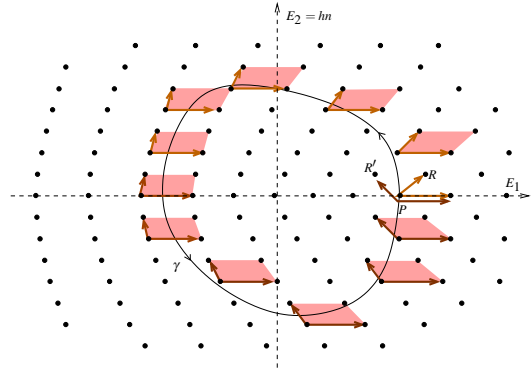

\subsection{Toric systems and polyads}

The simplest case of integrable systems whose global geometry is
perfectly understood is the toric one. Let $(f_1,\dots,f_n)$ be a
completely integrable system whose momentum map $F$ is proper.
\begin{defi}
\label{defi:torique}
  The system $F$ is of \textbf{toric type} if there exists an
  effective Hamiltonian action of $\T^n$ with momentum map $\Phi$ of
  the form $\Phi=\phy\circ F$ where $\phy$ is a local diffeomorphism
  on the image of $F$.
\end{defi}
Recall that by the Atiyah-Guillemin-Sternberg theorem
~\cite{atiyah-convex,guillemin-sternberg} the fibres of a momentum map
for a Hamiltonian torus action are connected and the image is a
rational convex polytope; and by Delzant's theorem~\cite{delzant} this
image fully determines the manifold $M$ and the momentum map $\Phi$
(up to isomorphism) in the completely integrable case. Actually the
connexity/convexity theorem is stated for compact $M$, but still holds
if the momentum map is proper~\cite{lerman-al} (the ``polytope'' being
not necessarily bounded). In this text we shall use abusively the
terms polytopes and polygons for possibly non-compact polyhedral
sets.

One can show that definition~\ref{defi:torique} entails that the
fibres of $F$ are connected and $\phy$ is a diffeomorphism from the
image of $F$ into the image of $\Phi$. Hence the structure of the
fibration is perfectly known. In particular the singularities are all
of (transversally) elliptic type. Thus by glueing local descriptions
one can easily obtain, in the semiclassical framework, a global
description of the joint spectrum.
\begin{theo}[\cite{san-redistribution}]
  \label{theo:BSglobal}
  Let $P=(P_1,\dots,P_n)$ be a quantum completely integrable system
  whose classical limit is of toric type, and let $\Sigma(P)$ be its
  joint spectrum. There exists a map $\phy_\h$ from $\RM^n$ to $\RM^n$
  such that for any compact $K\subset \RM^n$, 
  \begin{itemize}
  \item the restriction of $\phy_\h$ to $K$ is a classical symbol:
    $\phy_\h=\phy_0+\h\phy_1+\cdots$ whose principal term $\phy_0$ is
    a local diffeomorphism;
\item the components of $\phy_0\circ F$ are classical action
  variables;
\item  $\phy_{\h}(K\cap\Sigma(P))= \h\ZM^n\cap\phy_{\h}(K)+\ohb$.
  \end{itemize}
\end{theo}
The joint spectrum is thus transformed into a straight lattice
associated to the moment polytope, possibly shifted from the latter
because of the appearance of subprincipal terms.

Once a basis vector of the lattice in which sits the joint spectrum is
chosen, one can define a grouping of eigenvalues by associating those
that are on the same affine line directed by this vector. These
``packets'' of eigenvalues are called \textbf{polyads}, in reference
to the case of the harmonic oscillator in the physico-chemistry
literature. From the classical point of view, such a basis vector
defines a sub-action of $S^1$; the corresponding polyads are then
Weinstein's \emph{clusters}~\cite{weinstein-cluster}. Notice that in
our completely integrable situation, Weinstein's assumption on the
subprincipal symbols is no more necessary.

This theorem yields an amusing proof of the following corollary (which
can be proved by more usual means): when $M=T^*X$ (which is obviously
the case of the theorem, unless we state it in a Toeplitz context),
there can be only one critical point of maximal corank. This indicates
that in general the true interest of the theorem is actually not the
global description of the spectrum, which would hold for a very
limited set of examples, but instead the description, for a more
general joint spectrum, of all parts (convex, polyhedral) which
correspond to a sub-system of toric type.

\subsection{Almost toric systems}

A system of toric type is essentially a completely integrable system
all of whose singularities are non-degenerate and of (transversally)
elliptic type (although as it is, the assertion is not true: see
proposition ~\ref{prop:torique} below).

A mild way of generalising toric systems is to allow isolated
singularities in the momentum image.
\begin{defi}
  A completely integrable system with proper momentum map is
  \textbf{almost toric} when its singularities are all non-degenerate,
  without real hyperbolic block.
\end{defi}
In other words an almost toric system admits only elliptic or
focus-focus blocks. From now on we restrict to two degrees of freedom
systems (symplectic 4-manifolds). The classification up to
diffeomorphism of compact symplectic manifolds of dimension 4
admitting an almost toric system has just been carried out by Leung
and Symington~\cite{symington-four,leung-symington}.

Amongst elementary properties of almost toric systems, one can state:
\begin{prop}
    \label{prop:torique}
    \begin{itemize}
    \item If $F$ is a proper momentum map with non-degenerate
      singularities and the set of regular values of $F$ is connected,
      then $F$ is almost toric;
  \item if $F$ is almost toric then $F$ is of toric type if and only
    if the set of regular values of $F$ is connected and simply
    connected.
    \end{itemize}
\end{prop}

\paragraph{Generalised polytopes ---} At the time of writing this
article, no general result about semiclassics of almost toric systems
is known. To start with, let us consider a simple sub-class of two
degrees of freedom almost toric systems, namely those whose deficiency
index is 1:
\begin{defi}
  An almost toric integrable system $F=(f_1,f_2)$ on a symplectic
  manifold of dimension 4 has \textbf{deficiency index} equal to 1 if
  there exists a local diffeomorphism $\phy=(\phy_1,\phy_2)$ on the
  image of $F$ such that $\phy_1\circ F$ is a proper momentum map for
  an effective Hamiltonian action of $S^1$.
\end{defi}
Hamiltonian actions of $S^1$ on compact 4-manifolds have been
classified, both topologically~\cite{audin-topology} and
symplectically~\cite{karshon-S1}.  One can show that systems with
deficiency index 1 are often compact and hence subject to this
classification. However the symplectic manifold itself is not what
really matters here, considering that it is often obtained by a
``symplectic cutting'' adapted to the part of the system under study.
From the semiclassics viewpoint, the essential object is the image of
the momentum map, endowed with its structure of integral affine
manifold (with singularities).

Now let $F=(f_1,f_2)$ be an almost toric system with deficiency index
1; denote by $B\subset\RM^2$ the image of $F$, $B_r$ the set of
regular values, and $m_f$ the number of critical values
$c_1,\dots,c_{m_f}$ of focus-focus type. Let
$\vec\epsilon\in\{-1,+1\}^{m_f}$.  Denote by $\ell_i$ the vertical
half-line from $c_i$ in the direction given by $\epsilon_i$. Let $k_i$
be the \emph{monodromy index} of $c_i$ (it was shown
in~\cite{cushman-san} that in this situation the monodromy is
\emph{abelian} and can be identified to an integer valued index).
Let $\A^2_{\ZM}$ be the space $\RM^2$ equipped with the standard
integral affine structure.
\begin{theo}[\cite{san-polytope}]
There exists a homeomorphism $\psi$ from $B$ to $\psi(B)\subset\A^2_{\ZM}$
of the form   $\psi(x,y)=(x,\psi^{(2)}(x,y))$ such that
\begin{enumerate}
\item in the complement of $\cup_i\ell_i$, $\psi$ is an affine
  diffeomorphism (ie. the components of $\psi\circ F$ are local action
  variables)
\item $\psi$ extends to a $\Cinf$ multivalued map from $B_r$ to
  $\A^2_{\ZM}$ (with branching at all $\ell_i$) and for all
  $i=1,\dots,m_f$ and all $c\in\ell_i$,
\begin{equation*}
      \lim_{\substack{(x,y)\fleche c\\x<x_i}}d\psi(x,y) =
      \left(\begin{array}{cc}
          1 & 0\\ \epsilon_i k_i & 1
        \end{array}\right)
      \lim_{\substack{(x,y)\fleche
          c\\x>x_i}}d\psi(x,y),
      \label{equ:limit}
    \end{equation*}
  \item The image of $\psi$ is a rational convex polygon.
\end{enumerate}
\end{theo}
\begin{figure}[htbp]
  \begin{center}
    \input{cut.pstex_t}
    \caption{Construction of the polygon}
    \label{fig:cut}
  \end{center}
\end{figure}
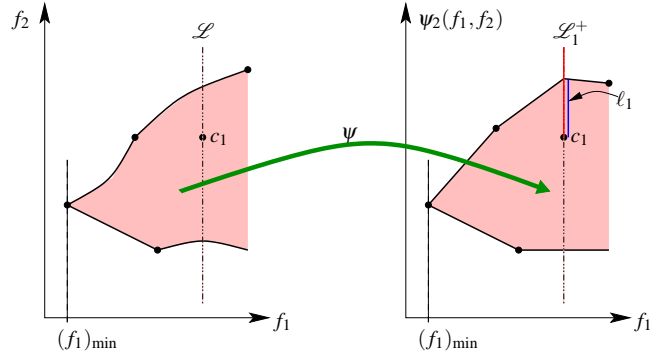
Although the system $F$ is not toric, the theorem still provides us
with a way of associating to it a rational convex polygon, which turns
out to be very useful to study the system. One can for instance write
localisation formulas that express Duistermaat-Heckman measures
associated with the action of the system in terms of the monodromy
index~\cite{san-polytope}.

The difference with the toric case, which of course is the main thrust
for these systems, is that the ``generalised'' moment polygon is not
unique; on the contrary, it is parameterised by a multi-sign
$\vec\epsilon$, which endows the class of possible polygons with an
abelian group structure and expresses the non-uniqueness of action
variables.

\subsection{Bifurcations and redistribution of eigenvalues}

Consider the following general question. Let an integrable quantum
system $P(t)$ depend on a parameter $t\in[0,1]$ such that $P(0)$ and
$P(1)$ are of toric type. We known how to describe the joint spectrum
of $P(0)$ and $P(1)$ thanks to theorem~\ref{theo:BSglobal}. What is
the relation between both spectra ?

To be more precise let us assume that there is way to define a
particular rational direction in the image of the joint principal
symbol $p(t)$, independently of $t$ (this is the case for instance of
there exists a sub-action of $S^1$ independent of $t$). One can then
define the corresponding \emph{polyads} for $P(0)$ and $P(1)$. How do
the eigenvalues rearrange from a set of polyad to another ? This is
the so-called \emph{redistribution problem}. According to
theorem~\ref{theo:BSglobal} it is enough to study the transformation
of the moment polygon to obtain the asymptotic behaviour of the number
of eigenvalues in each polyad.

By looking at the example of the coupling of two spins (it is a
Hamiltonian system on $S^2\times S^2$ which satisfies our hypothesis)
Sadovski\'{\i} and Zhilinski\'{\i}  conjectured that this
redistribution was related to the appearance of monodromy for certain
intermediate values of $t$~\cite{sadovski-zhilinski}.

Now assume that the system $p(t)$ is almost toric with deficiency
index 1, except for a finite number of $t$'s which we shall call
bifurcation times. Under the assumption that the only bifurcations
that the system undergoes are Hamiltonian Hopf
bifurcations\footnote{All these conditions are satisfied for the
  example of two spins.} (which correspond to a transformation
elliptic-elliptic $\leftrightarrow$ focus-focus and are
generic\footnote{They are generic for instance in the class of
  Hamiltonians that commute with a fixed $S^1$
  action}\cite{van-der-meer}), the conjecture is confirmed in the
following way: the polygons associated to $P(0)$ and $P(1)$ are
generalised polygons for a common system and the
$\vec\epsilon=(\epsilon_1,\dots,\epsilon_{m_f})$ corresponding to
their difference if determined by the sequence of Hopf bifurcations.
In other words one passes from a polygon to the other by a piecewise
affine transformation characterised by the position of bifurcating
critical values and their monodromy index.
\begin{figure}[htbp]
  \begin{center}
    \emph{Image of the momentum map:}
    
    \includegraphics[width=2.2cm]{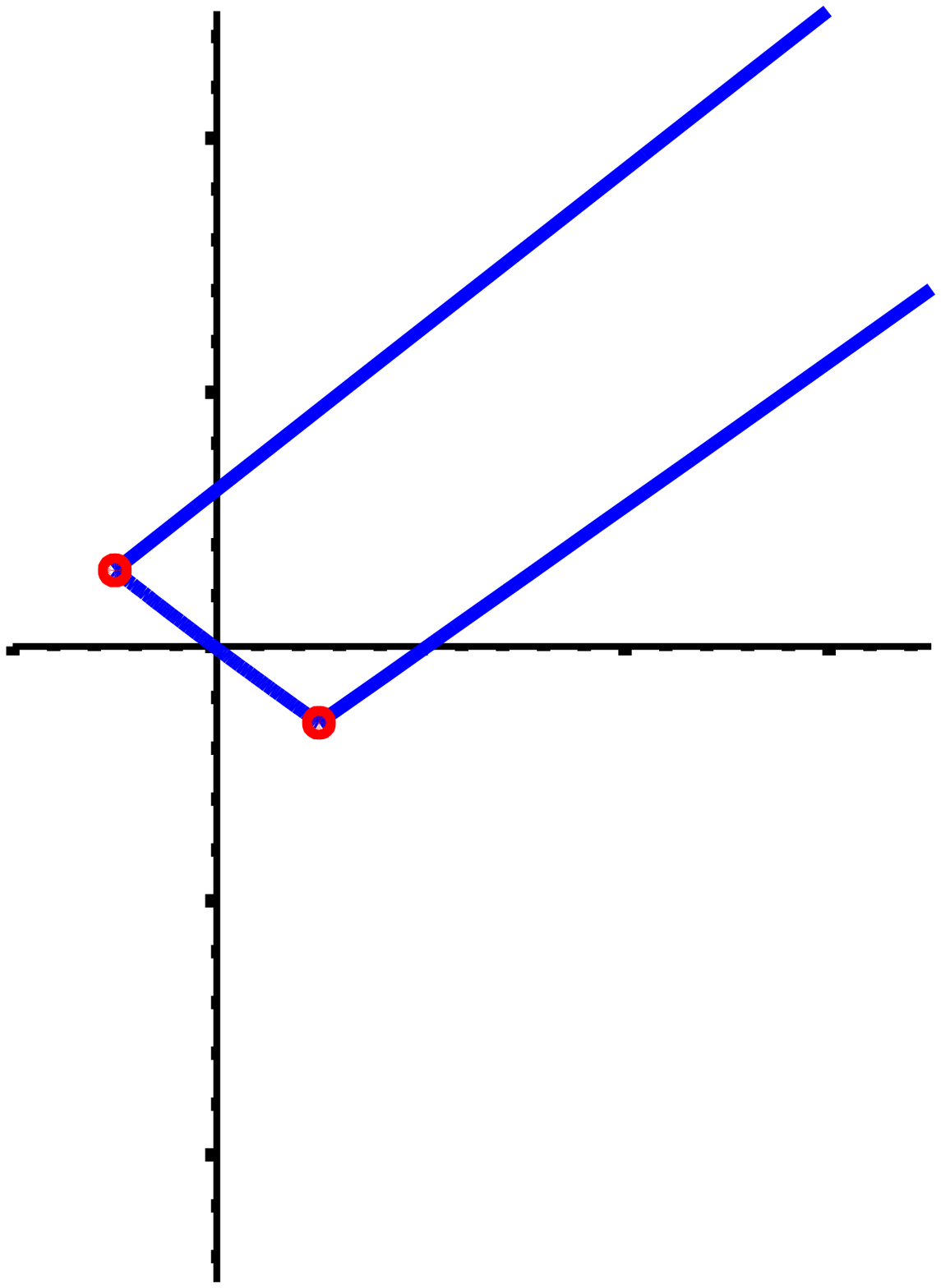}
    \includegraphics[width=2.2cm]{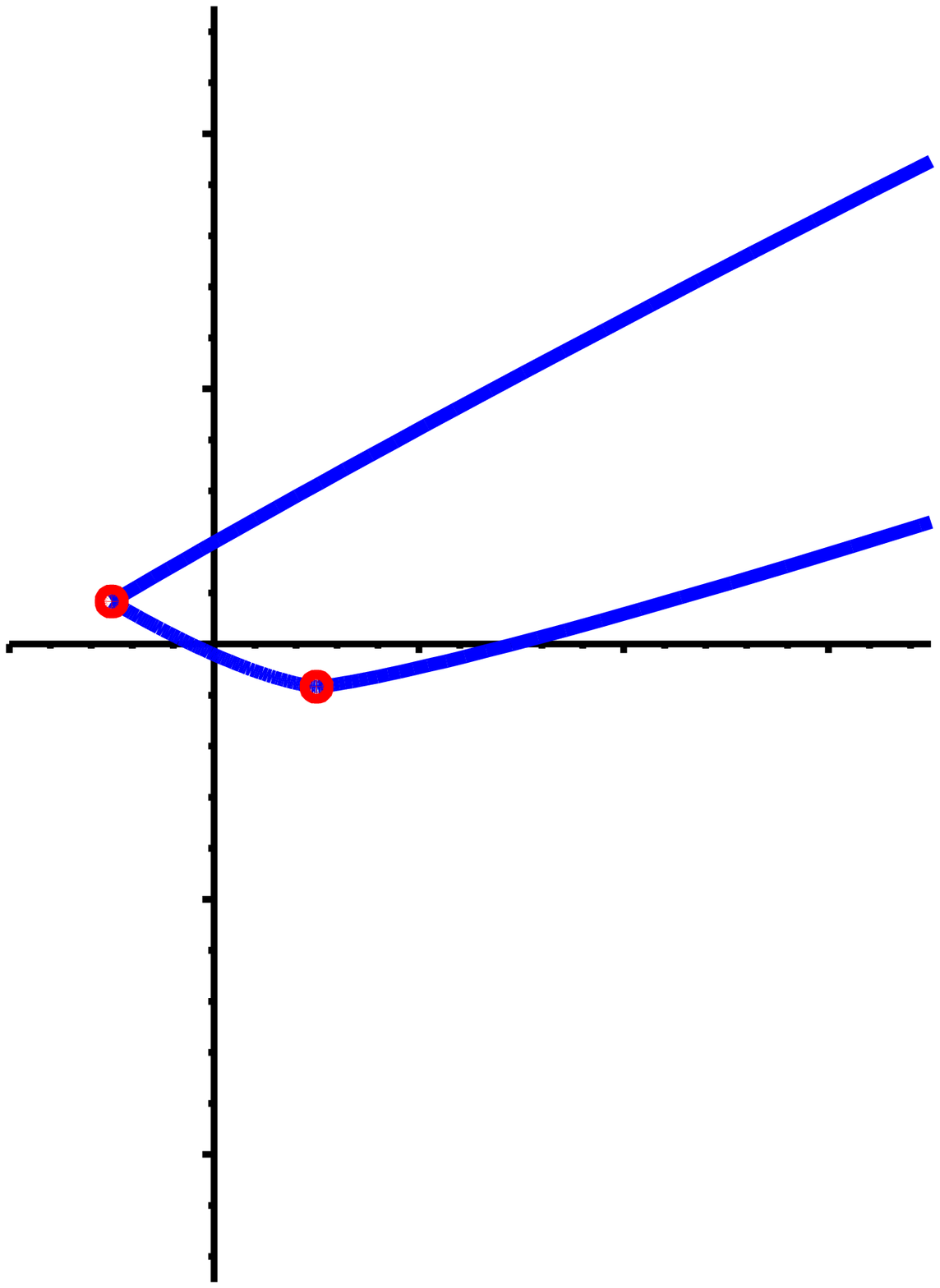}
    \includegraphics[width=2.2cm]{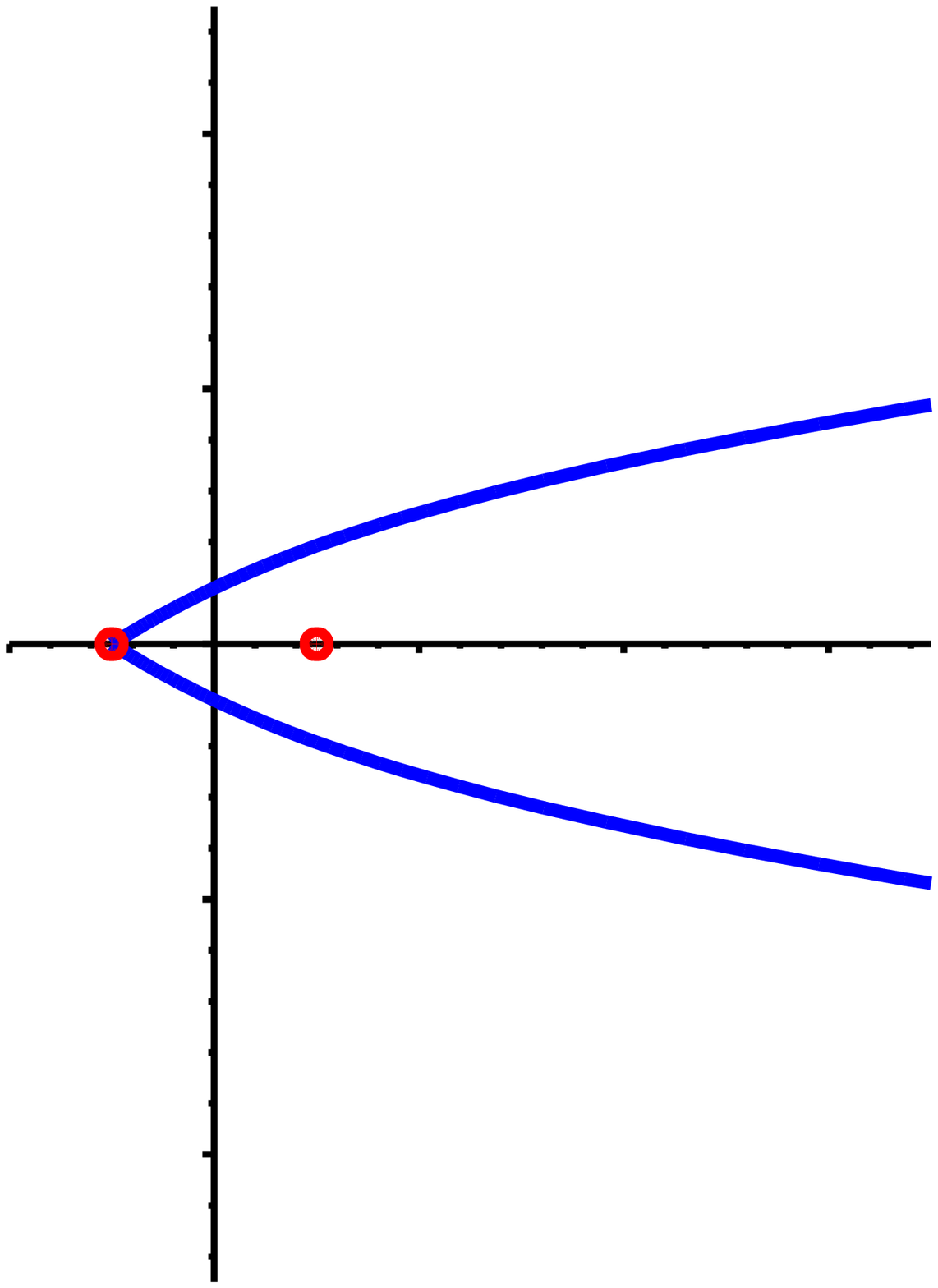}
    \includegraphics[width=2.2cm]{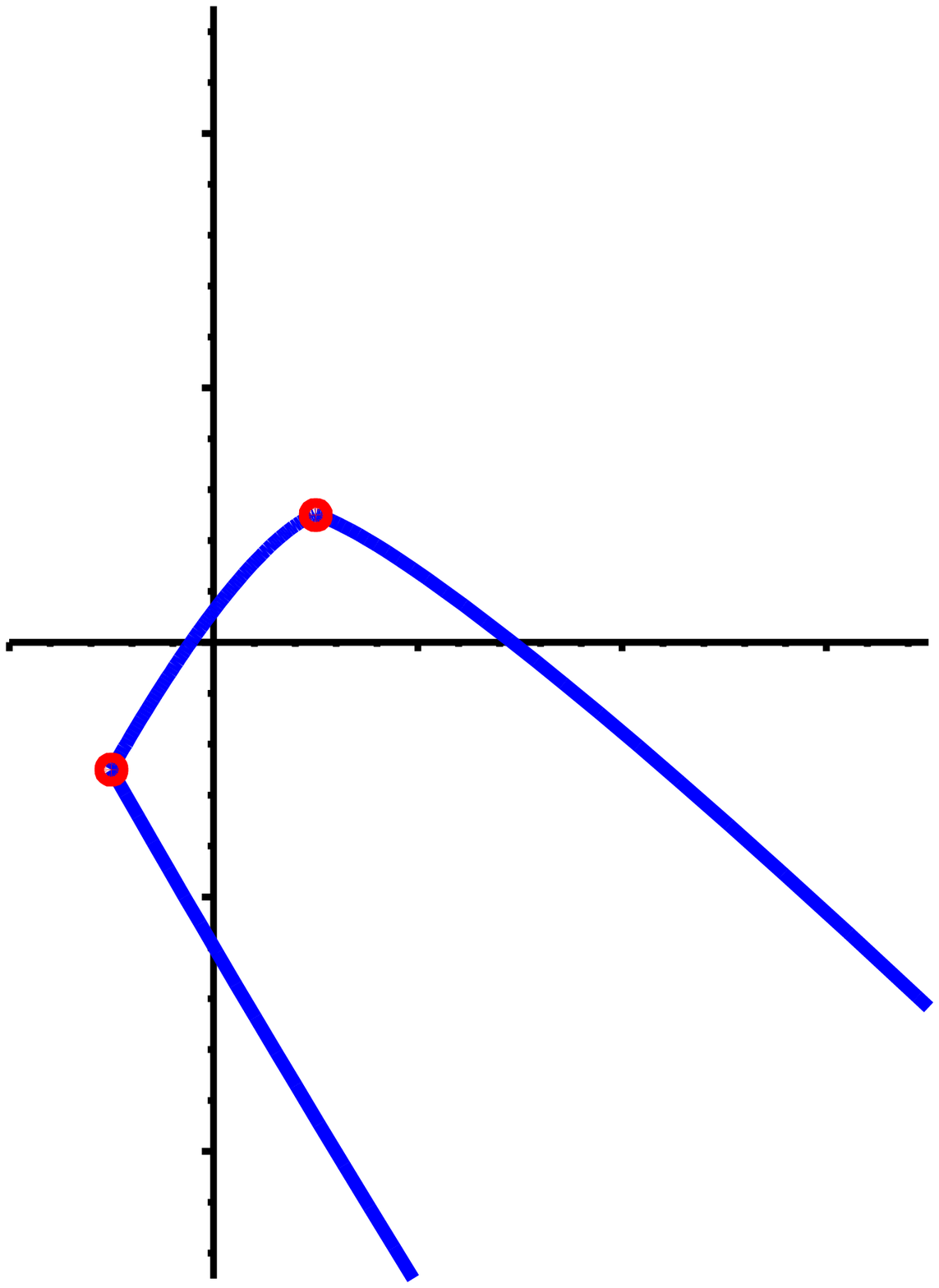}
    \includegraphics[width=2.2cm]{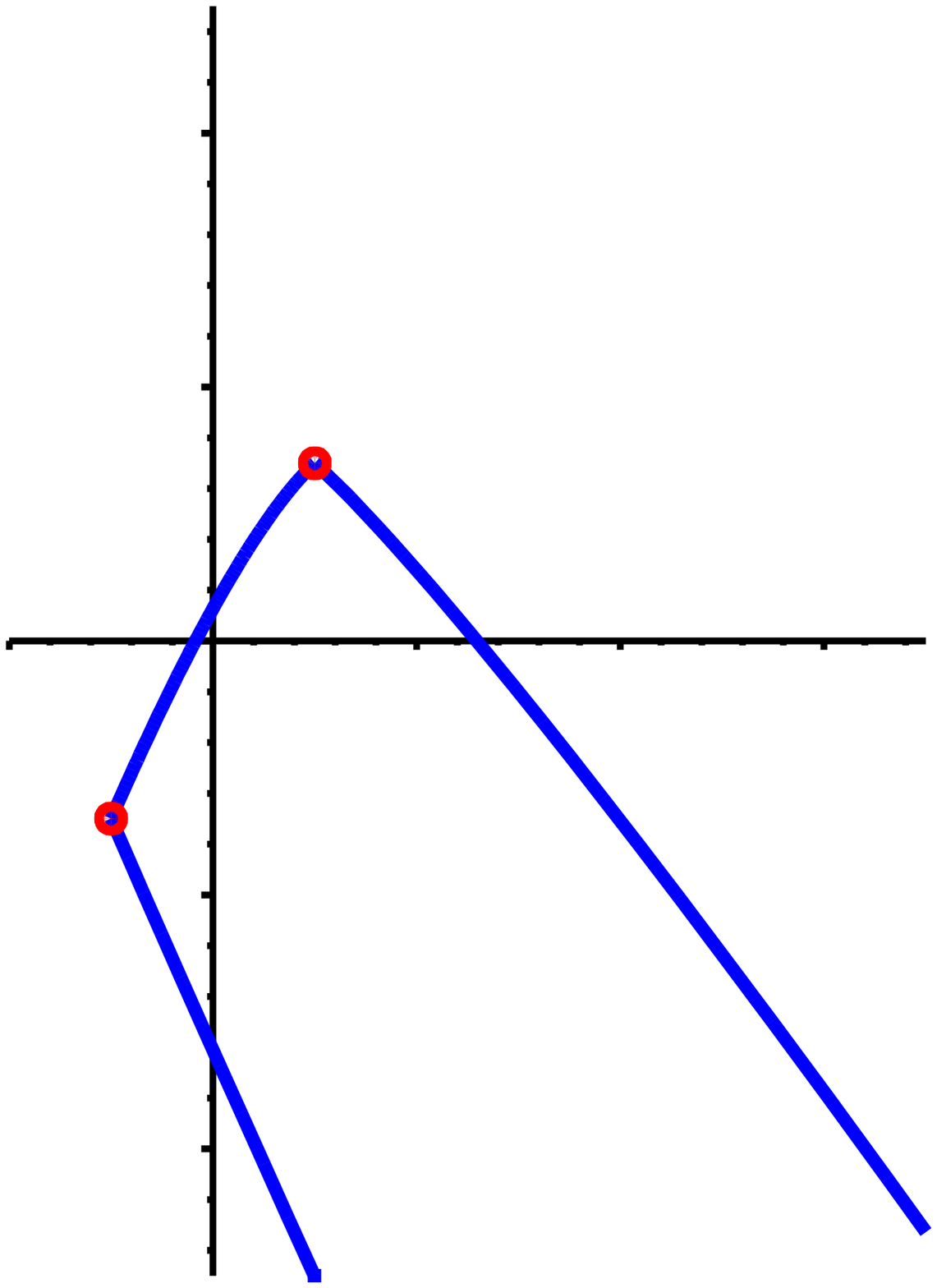}
    
    \emph{Corresponding generalised polytopes:}
    
    \includegraphics[width=2.2cm]{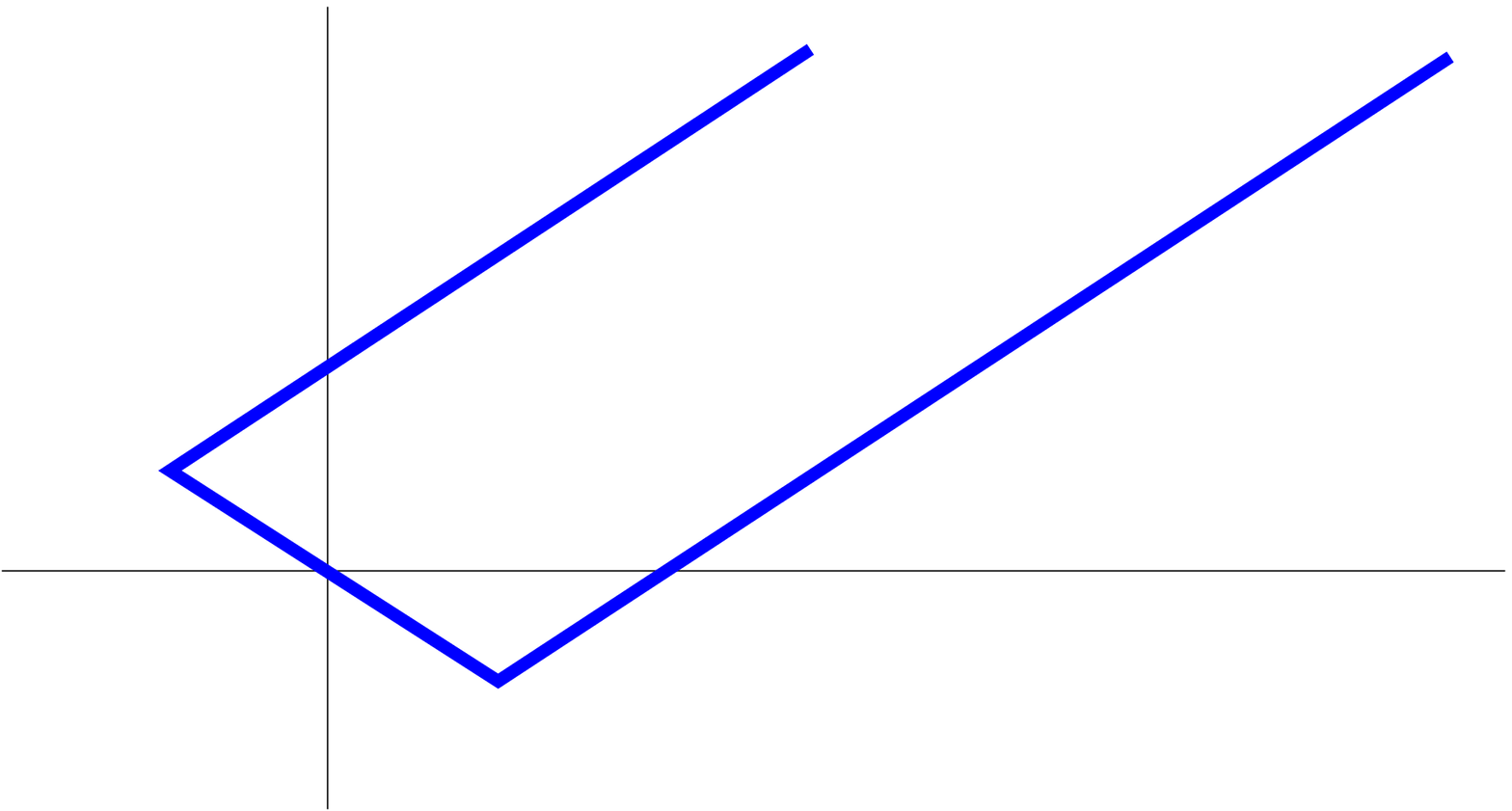}
    \includegraphics[width=2.2cm]{polytope+1}
    \fbox{\begin{minipage}[c]{2.2cm} \vspace{0pt}
        
        \includegraphics[width=\textwidth]{polytope+1}
      
        \centering\emph{and}
      
        \includegraphics[width=\textwidth]{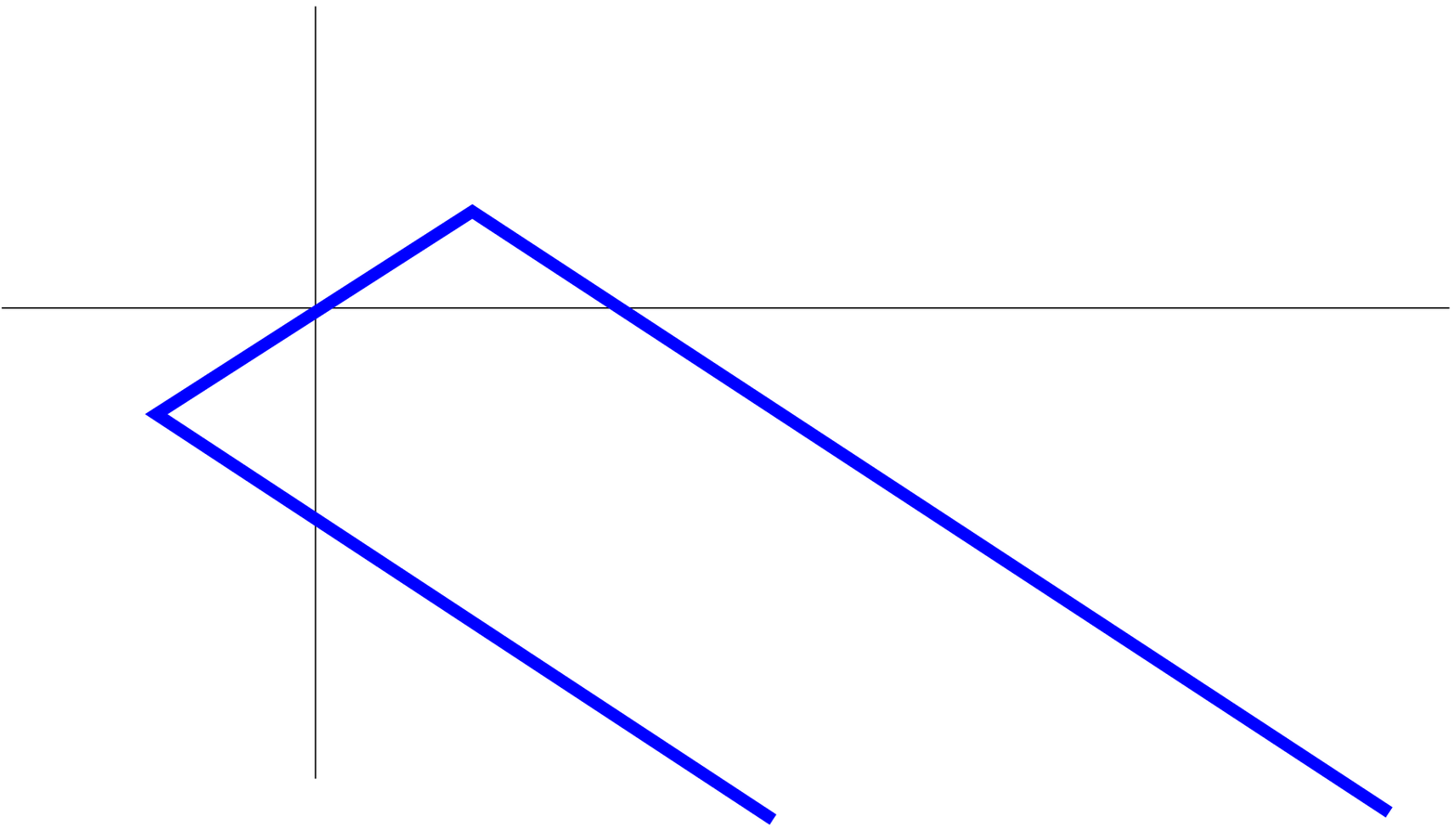}
    \end{minipage}}
  \includegraphics[width=2.2cm]{polytope-1}
  \includegraphics[width=2.2cm]{polytope-1}
    \caption{Bifurcation of the image of the momentum map for the
      coupling between a spin and a harmonic oscillator
      ($M=S^2\times\RM^2$)}
    \label{fig:moment_t}
  \end{center}
\end{figure}

This result was the main motivation for introducing these generalised
polytopes, since they perfectly describe in a geometrical and
combinatorial way way the eigenvalue redistribution amongst polyads.
It is natural however to imagine other applications of these
polytopes. Combined with the semi-global classification of
theorem~\ref{theo:invariants} they may turn out to be the right tool
for a global classification theorem \`a la Delzant.

\paragraph{Acknowledgements ---} This article was written during a
very pleasant stay at the Centre Bernoulli (EPFL, Lausanne). I am very
grateful to Tudor Ratiu for his kind invitation. This work was
partially supported by the european network MASIE.

\setlength{\parskip}{0ex}
\bibliographystyle{plain}%
\bibliography{bibli}

\end{document}

%% file: huit.pstex_t
\begin{picture}(0,0)%
\includegraphics{huit.pstex}%
\end{picture}%
\setlength{\unitlength}{3315sp}%
\begingroup\makeatletter\ifx\SetFigFont\undefined%
\gdef\SetFigFont#1#2#3#4#5{%
  \reset@font\fontsize{#1}{#2pt}%
  \fontfamily{#3}\fontseries{#4}\fontshape{#5}%
  \selectfont}%
\fi\endgroup%
\begin{picture}(5299,2001)(484,-2005)
\put(3106,-1951){\makebox(0,0)[lb]{\smash{\SetFigFont{10}{12.0}{\rmdefault}{\mddefault}{\updefault}\special{ps: gsave 0 0 0 setrgbcolor}topological\special{ps: grestore}}}}
\put(4366,-1951){\makebox(0,0)[lb]{\smash{\SetFigFont{10}{12.0}{\rmdefault}{\mddefault}{\updefault}\special{ps: gsave 0 0 0 setrgbcolor}analytic\special{ps: grestore}}}}
\put(5491,-1951){\makebox(0,0)[lb]{\smash{\SetFigFont{10}{12.0}{\rmdefault}{\mddefault}{\updefault}\special{ps: gsave 0 0 0 setrgbcolor}$\Cinf$\special{ps: grestore}}}}
\put(3916,-241){\makebox(0,0)[lb]{\smash{\SetFigFont{10}{12.0}{\rmdefault}{\mddefault}{\updefault}\special{ps: gsave 0 0 0 setrgbcolor}Leaf space\special{ps: grestore}}}}
\end{picture}

%% file: jspec-hyp.pstex_t
\begin{picture}(0,0)%
\includegraphics{jspec-hyp.pstex}%
\end{picture}%
\setlength{\unitlength}{2486sp}%
\begingroup\makeatletter\ifx\SetFigFont\undefined%
\gdef\SetFigFont#1#2#3#4#5{%
  \reset@font\fontsize{#1}{#2pt}%
  \fontfamily{#3}\fontseries{#4}\fontshape{#5}%
  \selectfont}%
\fi\endgroup%
\begin{picture}(3722,3831)(644,-3751)
\put(1891,-3751){\makebox(0,0)[lb]{\smash{\SetFigFont{8}{9.6}{\rmdefault}{\mddefault}{\updefault}\special{ps: gsave 0 0 0 setrgbcolor}$\h$\special{ps: grestore}}}}
\put(4366,-1366){\makebox(0,0)[lb]{\smash{\SetFigFont{8}{9.6}{\rmdefault}{\mddefault}{\updefault}\special{ps: gsave 0 0 0 setrgbcolor}$\h/\ln\h$\special{ps: grestore}}}}
\end{picture}

%% file: monodromie.pstex_t
\begin{picture}(0,0)%
\includegraphics{monodromie.pstex}%
\end{picture}%
\setlength{\unitlength}{1776sp}%
\begingroup\makeatletter\ifx\SetFigFont\undefined%
\gdef\SetFigFont#1#2#3#4#5{%
  \reset@font\fontsize{#1}{#2pt}%
  \fontfamily{#3}\fontseries{#4}\fontshape{#5}%
  \selectfont}%
\fi\endgroup%
\begin{picture}(7271,5124)(169,-6073)
\put(4051,-1186){\makebox(0,0)[lb]{\smash{\SetFigFont{5}{6.0}{\rmdefault}{\mddefault}{\updefault}\special{ps: gsave 0 0 0 setrgbcolor}$E_2=hn$\special{ps: grestore}}}}
\put(5612,-3898){\makebox(0,0)[lb]{\smash{\SetFigFont{5}{6.0}{\rmdefault}{\mddefault}{\updefault}\special{ps: gsave 0 0 0 setrgbcolor}$P$\special{ps: grestore}}}}
\put(6032,-3411){\makebox(0,0)[lb]{\smash{\SetFigFont{5}{6.0}{\rmdefault}{\mddefault}{\updefault}\special{ps: gsave 0 0 0 setrgbcolor}$R$\special{ps: grestore}}}}
\put(5140,-3306){\makebox(0,0)[lb]{\smash{\SetFigFont{5}{6.0}{\rmdefault}{\mddefault}{\updefault}\special{ps: gsave 0 0 0 setrgbcolor}$R'$\special{ps: grestore}}}}
\put(2461,-4664){\makebox(0,0)[lb]{\smash{\SetFigFont{5}{6.0}{\rmdefault}{\mddefault}{\updefault}\special{ps: gsave 0 0 0 setrgbcolor}$\gamma$\special{ps: grestore}}}}
\put(7111,-3526){\makebox(0,0)[lb]{\smash{\SetFigFont{5}{6.0}{\rmdefault}{\mddefault}{\updefault}\special{ps: gsave 0 0 0 setrgbcolor}$E_1$\special{ps: grestore}}}}
\end{picture}

%% file: cut.pstex_t
\begin{picture}(0,0)%
\includegraphics{cut.pstex}%
\end{picture}%
\setlength{\unitlength}{2486sp}%
\begingroup\makeatletter\ifx\SetFigFont\undefined%
\gdef\SetFigFont#1#2#3#4#5{%
  \reset@font\fontsize{#1}{#2pt}%
  \fontfamily{#3}\fontseries{#4}\fontshape{#5}%
  \selectfont}%
\fi\endgroup%
\begin{picture}(6210,3490)(361,-2864)
\put(4411,-2806){\makebox(0,0)[lb]{\smash{\SetFigFont{8}{9.6}{\rmdefault}{\mddefault}{\updefault}\special{ps: gsave 0 0 0 setrgbcolor}$(f_1)_{\textup{min}}$\special{ps: grestore}}}}
\put(5941,-781){\makebox(0,0)[lb]{\smash{\SetFigFont{8}{9.6}{\rmdefault}{\mddefault}{\updefault}\special{ps: gsave 0 0 0 setrgbcolor}$c_1$\special{ps: grestore}}}}
\put(811,-2806){\makebox(0,0)[lb]{\smash{\SetFigFont{8}{9.6}{\rmdefault}{\mddefault}{\updefault}\special{ps: gsave 0 0 0 setrgbcolor}$(f_1)_{\textup{min}}$\special{ps: grestore}}}}
\put(2341,-781){\makebox(0,0)[lb]{\smash{\SetFigFont{8}{9.6}{\rmdefault}{\mddefault}{\updefault}\special{ps: gsave 0 0 0 setrgbcolor}$c_1$\special{ps: grestore}}}}
\put(2971,-2626){\makebox(0,0)[lb]{\smash{\SetFigFont{8}{9.6}{\rmdefault}{\mddefault}{\updefault}\special{ps: gsave 0 0 0 setrgbcolor}$f_1$\special{ps: grestore}}}}
\put(361,389){\makebox(0,0)[lb]{\smash{\SetFigFont{8}{9.6}{\rmdefault}{\mddefault}{\updefault}\special{ps: gsave 0 0 0 setrgbcolor}$f_2$\special{ps: grestore}}}}
\put(5761,254){\makebox(0,0)[lb]{\smash{\SetFigFont{8}{9.6}{\rmdefault}{\mddefault}{\updefault}\special{ps: gsave 0 0 0 setrgbcolor}$\mathcal{L}_1^+$\special{ps: grestore}}}}
\put(2161,254){\makebox(0,0)[lb]{\smash{\SetFigFont{8}{9.6}{\rmdefault}{\mddefault}{\updefault}\special{ps: gsave 0 0 0 setrgbcolor}$\mathcal{L}$\special{ps: grestore}}}}
\put(3646,-736){\makebox(0,0)[lb]{\smash{\SetFigFont{8}{9.6}{\rmdefault}{\mddefault}{\updefault}\special{ps: gsave 0 0 0 setrgbcolor}$\psi$\special{ps: grestore}}}}
\put(6391,-421){\makebox(0,0)[lb]{\smash{\SetFigFont{8}{9.6}{\rmdefault}{\mddefault}{\updefault}\special{ps: gsave 0 0 0 setrgbcolor}$\ell_1$\special{ps: grestore}}}}
\put(6571,-2626){\makebox(0,0)[lb]{\smash{\SetFigFont{8}{9.6}{\rmdefault}{\mddefault}{\updefault}\special{ps: gsave 0 0 0 setrgbcolor}$f_1$\special{ps: grestore}}}}
\put(4411,389){\makebox(0,0)[lb]{\smash{\SetFigFont{8}{9.6}{\rmdefault}{\mddefault}{\updefault}\special{ps: gsave 0 0 0 setrgbcolor}$\psi_2(f_1,f_2)$\special{ps: grestore}}}}
\end{picture}